\documentclass[a4paper]{amsart}
\usepackage{amssymb,amsmath,amsthm,mathrsfs,fancyhdr,lastpage}
\usepackage[dvips]{color}

\def\red#1{\textcolor{red}{#1}}

\newcommand{\mythicklines}{\linethickness {0.40mm}}
\newcommand{\mythinlines}{\linethickness {0.12mm}}

\newtheorem{theorem}{Theorem}[section]
\newtheorem{lemma}[theorem]{Lemma}
\newtheorem{proposition}[theorem]{Proposition}
\newtheorem*{sublem}{Sublemma}
\newtheorem*{pfs}{\rm {\it Proof of Sublemma}}
\numberwithin{equation}{section}

\begin{document}

\title[Arc index of pretzel knots of type $(-p, q, r)$]{Arc index of pretzel knots of type \boldmath$(-p, q, r)$}
\author[H. J. Lee]{Hwa Jeong Lee}
\address{Department of Mathematical Sciences, KAIST, 291 Daehak-ro,
Yuseong-gu, Daejeon 305-701, Korea}
\email{hjwith@kaist.ac.kr}
\thanks{The first author was supported in part by the National Research Foundation of Korea Grant
funded by the Korean Government (NRF-2010-0024630)}

\author[G. T. Jin]{Gyo Taek Jin}
\address{Department of Mathematical Sciences, KAIST, 291 Daehak-ro,
Yuseong-gu, Daejeon 305-701, Korea}
\email{trefoil@kaist.ac.kr}
\thanks{The second author was supported in part by the National Research Foundation of Korea Grant
funded by the Korean Government (NRF-2011-0027989)}

\keywords{knot, pretzel knot, arc presentation, arc index, Kauffman polynomial}
\subjclass[2000]{Primary 57M27; Secondary 57M25}

\maketitle

\begin{abstract}
We computed the arc index for some of the pretzel knots $K=P(-p,q,r)$ with $p,q,r\ge2$, $r\geq q$ and at most one of $p,q,r$ is even. If $q=2$, then the arc index $\alpha(K)$ equals the minimal crossing number $c(K)$. If $p\ge3$ and $q=3$, then $\alpha(K)=c(K)-1$. If $p\ge5$ and $q=4$, then $\alpha(K)=c(K)-2$.
\end{abstract}


\section{Arc presentation}
\begin{figure}[ht]
\centering
\begin{picture}(60,60)
\put(0,0){\line(1,0){40}}
\put(0,0){\line(0,1){40}}
\put(40,0){\line(0,1){17}} \put(40,23){\line(0,1){37}}
\put(40,60){\line(-1,0){20}}
\put(20,60){\line(0,-1){17}} \put(20,37){\line(1,-1){17}}
\put(37,20){\line(1,0){23}}
\put(60,20){\line(0,1){20}}
\put(60,40){\line(-1,0){17}} \put(37,40){\line(-1,0){37}}
\mythicklines
\qbezier(28.5,28.5)(50,28.5)(50,40)
\qbezier(28.5,28.5)(28.5,50)(40,50)
\qbezier(50,40)(50,50)(40,50)
\end{picture}
\caption{An arc presentation of the right-handed trefoil knot}\label{fig:arc pres of trefoil}
\end{figure}
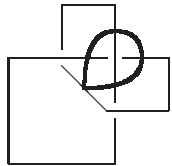

Let $D$ be a diagram of a knot or a link $L$. Suppose that there is a simple closed curve $C$ meeting $D$ in $k$ distinct points which divide $D$ into $k$ arcs $\alpha_1,\alpha_2,\ldots,\alpha_k$ with the following properties:
\begin{enumerate}
\item Each $\alpha_i$ has no self-crossing.
\item If $\alpha_i$ crosses over $\alpha_j$ at a crossing, then $i>j$ and it crosses over $\alpha_j$ at any other crossings with $\alpha_j$.
\item For each $i$, there exists an embedded disk $d_i$ such that $\partial d_i=C$ and $\alpha_i\subset d_i$.
\item $d_i\cap d_j=C$, for distinct $i$ and $j$.
\end{enumerate}
Then the pair $(D,C)$ is called an \emph{arc presentation\/} of $L$ with $k$ arcs, and $C$ is called the \emph{axis\/} of the arc presentation.
Figure~\ref{fig:arc pres of trefoil} shows an arc presentation of the trefoil knot. The thick round curve is the axis. It is known that every knot or link has an arc presentation~\cite{Brunn,C1995}. For a given knot or link $L$, the minimal number of arcs in all arc presentations of $L$ is called the \emph{arc index\/} of $L$, denoted by $\alpha(L)$.

By removing a point $P$ from $C$ away from $L$, we may identify $C\setminus P$ with the $z$-axis and each $d_i\setminus P$ with a vertical half plane along the $z$-axis. This shows that an arc presentation is equivalent to an open-book presentation.
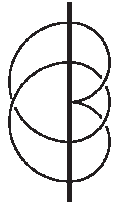
\begin{figure}[ht]
\setlength{\unitlength}{1.5pt}
\centering
\begin{picture}(50,50)(-25,-25)
\qbezier(0.,0.)(0.6554286300,0.000091492)(1.305261921,0.085551386)
\qbezier(1.305261921,0.085551386)(1.955071329,0.171192699)(2.588190451,0.340741737)
\qbezier(2.588190451,0.340741737)(3.221262214,0.510467524)(3.826834325,0.761204675)
\qbezier(3.826834325,0.761204675)(4.432336411,1.012110881)(5.,1.339745960)
\qbezier(5.,1.339745960)(5.567572098,1.667539511)(6.087614288,2.066466596)
\qbezier(6.087614288,2.066466596)(6.607545088,2.465538852)(7.071067810,2.928932190)
\qbezier(7.071067810,2.928932190)(7.534461147,3.392454911)(7.933533403,3.912385710)
\qbezier(7.933533403,3.912385710)(8.332460488,4.432427900)(8.660254040,5.)
\qbezier(8.660254040,5.)(8.987889120,5.567663591)(9.238795325,6.173165675)
\qbezier(9.238795325,6.173165675)(9.489532475,6.778737785)(9.659258263,7.411809549)
\qbezier(9.659258263,7.411809549)(9.828807301,8.044928670)(9.914448614,8.694738078)
\qbezier(9.914448614,8.694738078)(9.999908508,9.344571368)(10.,10.)
\qbezier(10.,10.)(9.999908508,10.65542863)(9.914448614,11.30526192)
\qbezier(9.914448614,11.30526192)(9.828807301,11.95507133)(9.659258263,12.58819045)
\qbezier(9.659258263,12.58819045)(9.489532475,13.22126222)(9.238795325,13.82683432)
\qbezier(9.238795325,13.82683432)(8.987889120,14.43233641)(8.660254040,15.)
\qbezier(8.660254040,15.)(8.332460488,15.56757210)(7.933533403,16.08761429)
\qbezier(7.933533403,16.08761429)(7.534461147,16.60754509)(7.071067810,17.07106781)
\qbezier(7.071067810,17.07106781)(6.607545088,17.53446115)(6.087614288,17.93353340)
\qbezier(6.087614288,17.93353340)(5.567572098,18.33246049)(5.,18.66025404)
\qbezier(5.,18.66025404)(4.432336411,18.98788912)(3.826834325,19.23879532)
\qbezier(3.826834325,19.23879532)(3.221262214,19.48953248)(2.588190451,19.65925826)
\qbezier(2.588190451,19.65925826)(1.955071329,19.82880730)(1.305261921,19.91444861)
\qbezier(1.305261921,19.91444861)(0.6554286300,19.99990851)(0.,20.)
\qbezier(0.,-10.)(0.6554286300,-9.999908508)(1.305261921,-9.914448614)
\qbezier(1.305261921,-9.914448614)(1.955071329,-9.828807301)(2.588190451,-9.659258263)
\qbezier(2.588190451,-9.659258263)(3.221262214,-9.489532476)(3.826834325,-9.238795325)
\qbezier(3.826834325,-9.238795325)(4.432336411,-8.987889119)(5.,-8.660254040)
\qbezier(5.,-8.660254040)(5.567572098,-8.332460489)(6.087614288,-7.933533404)
\qbezier(6.087614288,-7.933533404)(6.607545088,-7.534461148)(7.071067810,-7.071067810)
\qbezier(7.071067810,-7.071067810)(7.534461147,-6.607545089)(7.933533403,-6.087614290)
\qbezier(7.933533403,-6.087614290)(8.332460488,-5.567572100)(8.660254040,-5.)
\qbezier(8.660254040,-5.)(8.987889120,-4.432336409)(9.238795325,-3.826834325)
\qbezier(9.238795325,-3.826834325)(9.489532475,-3.221262215)(9.659258263,-2.588190451)
\qbezier(9.659258263,-2.588190451)(9.828807301,-1.955071330)(9.914448614,-1.305261922)
\qbezier(9.914448614,-1.305261922)(9.999908508,-0.6554286315)(10.,0.)
\qbezier(10.,0.)(9.999908508,0.6554286315)(9.914448614,1.305261922)
\qbezier(9.914448614,1.305261922)(9.828807301,1.955071330)(9.659258263,2.588190451)
\qbezier(9.659258263,2.588190451)(9.489532475,3.221262215)(9.238795325,3.826834325)
%
\qbezier(7.933533403,6.087614290)(7.534461147,6.607545089)(7.071067810,7.071067810)
\qbezier(7.071067810,7.071067810)(6.607545088,7.534461148)(6.087614288,7.933533404)
\qbezier(6.087614288,7.933533404)(5.567572098,8.332460489)(5.,8.660254040)
\qbezier(5.,8.660254040)(4.432336411,8.987889119)(3.826834325,9.238795325)
\qbezier(3.826834325,9.238795325)(3.221262214,9.489532476)(2.588190451,9.659258263)
\qbezier(2.588190451,9.659258263)(1.955071329,9.828807301)(1.305261921,9.914448614)
\qbezier(1.305261921,9.914448614)(0.6554286300,9.999908508)(0.,10.)
\qbezier(0.,-20.)(0.6554286300,-19.99990851)(1.305261921,-19.91444861)
\qbezier(1.305261921,-19.91444861)(1.955071329,-19.82880730)(2.588190451,-19.65925826)
\qbezier(2.588190451,-19.65925826)(3.221262214,-19.48953248)(3.826834325,-19.23879532)
\qbezier(3.826834325,-19.23879532)(4.432336411,-18.98788912)(5.,-18.66025404)
\qbezier(5.,-18.66025404)(5.567572098,-18.33246049)(6.087614288,-17.93353340)
\qbezier(6.087614288,-17.93353340)(6.607545088,-17.53446115)(7.071067810,-17.07106781)
\qbezier(7.071067810,-17.07106781)(7.534461147,-16.60754509)(7.933533403,-16.08761429)
\qbezier(7.933533403,-16.08761429)(8.332460488,-15.56757210)(8.660254040,-15.)
\qbezier(8.660254040,-15.)(8.987889120,-14.43233641)(9.238795325,-13.82683432)
\qbezier(9.238795325,-13.82683432)(9.489532475,-13.22126222)(9.659258263,-12.58819045)
\qbezier(9.659258263,-12.58819045)(9.828807301,-11.95507133)(9.914448614,-11.30526192)
\qbezier(9.914448614,-11.30526192)(9.999908508,-10.65542863)(10.,-10.)
\qbezier(10.,-10.)(9.999908508,-9.344571368)(9.914448614,-8.694738078)
\qbezier(9.914448614,-8.694738078)(9.828807301,-8.044928670)(9.659258263,-7.411809549)
\qbezier(9.659258263,-7.411809549)(9.489532475,-6.778737785)(9.238795325,-6.173165675)
\qbezier(7.933533403,-3.912385710)(7.534461147,-3.392454911)(7.071067810,-2.928932190)
\qbezier(7.071067810,-2.928932190)(6.607545088,-2.465538852)(6.087614288,-2.066466596)
\qbezier(6.087614288,-2.066466596)(5.567572098,-1.667539511)(5.,-1.339745960)
\qbezier(5.,-1.339745960)(4.432336411,-1.012110881)(3.826834325,-0.761204675)
\qbezier(3.826834325,-0.761204675)(3.221262214,-0.510467524)(2.588190451,-0.340741737)
\qbezier(2.588190451,-0.340741737)(1.955071329,-0.171192699)(1.305261921,-0.085551386)
\qbezier(1.305261921,-0.085551386)(0.6554286300,-0.000091492)(0.,0.)
\qbezier(0.,20.)(-0.4910488796,19.99999140)(-0.9810469362,19.96788385)
\qbezier(-0.9810469362,19.96788385)(-1.471043868,19.93575914)(-1.957892882,19.87167292)
\qbezier(-1.957892882,19.87167292)(-2.444739636,19.80756966)(-2.926354830,19.71177921)
\qbezier(-2.926354830,19.71177921)(-3.407966654,19.61597189)(-3.882285676,19.48888739)
\qbezier(-3.882285676,19.48888739)(-4.356600236,19.36178629)(-4.821591980,19.20395194)
\qbezier(-4.821591980,19.20395194)(-5.286578186,19.04610132)(-5.740251488,18.85819299)
\qbezier(-5.740251488,18.85819299)(-6.193918197,18.67026878)(-6.634330356,18.45309112)
\qbezier(-6.634330356,18.45309112)(-7.074734898,18.23589805)(-7.500000000,17.99038106)
\qbezier(-7.500000000,17.99038106)(-7.925256498,17.74484918)(-8.333553494,17.47204419)
\qbezier(-8.333553494,17.47204419)(-8.741840937,17.19922490)(-9.131421432,16.90030011)
\qbezier(-9.131421432,16.90030011)(-9.520991466,16.60136167)(-9.890187226,16.27759711)
\qbezier(-9.890187226,16.27759711)(-10.25937165,15.95381963)(-10.60660172,15.60660172)
\qbezier(-10.60660172,15.60660172)(-10.95381963,15.25937165)(-11.27759711,14.89018723)
\qbezier(-11.27759711,14.89018723)(-11.60136167,14.52099147)(-11.90030010,14.13142144)
\qbezier(-11.90030010,14.13142144)(-12.19922490,13.74184094)(-12.47204418,13.33355350)
\qbezier(-12.47204418,13.33355350)(-12.74484918,12.92525650)(-12.99038106,12.50000000)
\qbezier(-12.99038106,12.50000000)(-13.23589805,12.07473490)(-13.45309112,11.63433035)
\qbezier(-13.45309112,11.63433035)(-13.67026878,11.19391820)(-13.85819299,10.74025149)
\qbezier(-13.85819299,10.74025149)(-14.04610132,10.28657819)(-14.20395194,9.821591981)
\qbezier(-14.20395194,9.821591981)(-14.36178629,9.356600243)(-14.48888739,8.882285676)
\qbezier(-14.48888739,8.882285676)(-14.61597189,8.407966662)(-14.71177921,7.926354832)
\qbezier(-14.71177921,7.926354832)(-14.80756966,7.444739645)(-14.87167292,6.957892883)
\qbezier(-14.87167292,6.957892883)(-14.93575914,6.471043878)(-14.96788385,5.981046938)
\qbezier(-14.96788385,5.981046938)(-14.99999140,5.491048875)(-15.,5.)
\qbezier(-15.,5.)(-14.99999140,4.508951125)(-14.96788385,4.018953062)
\qbezier(-14.96788385,4.018953062)(-14.93575914,3.528956122)(-14.87167292,3.042107117)
\qbezier(-14.87167292,3.042107117)(-14.80756966,2.555260355)(-14.71177921,2.073645168)
\qbezier(-14.71177921,2.073645168)(-14.61597189,1.592033338)(-14.48888739,1.117714324)
\qbezier(-13.45309112,-1.634330353)(-13.23589805,-2.074734904)(-12.99038106,-2.500000000)
\qbezier(-12.99038106,-2.500000000)(-12.74484918,-2.925256501)(-12.47204418,-3.333553496)
\qbezier(-12.47204418,-3.333553496)(-12.19922490,-3.741840940)(-11.90030010,-4.131421435)
\qbezier(-11.90030010,-4.131421435)(-11.60136167,-4.520991469)(-11.27759711,-4.890187228)
\qbezier(-11.27759711,-4.890187228)(-10.95381963,-5.25937165)(-10.60660172,-5.60660172)
\qbezier(-10.60660172,-5.60660172)(-10.25937165,-5.95381963)(-9.890187226,-6.27759711)
\qbezier(-9.890187226,-6.27759711)(-9.520991466,-6.60136167)(-9.131421432,-6.90030011)
\qbezier(-9.131421432,-6.90030011)(-8.741840937,-7.19922490)(-8.333553494,-7.47204419)
\qbezier(-8.333553494,-7.47204419)(-7.925256498,-7.74484918)(-7.500000000,-7.99038106)
\qbezier(-7.500000000,-7.99038106)(-7.074734898,-8.23589805)(-6.634330356,-8.45309112)
\qbezier(-6.634330356,-8.45309112)(-6.193918197,-8.67026878)(-5.740251488,-8.85819299)
\qbezier(-5.740251488,-8.85819299)(-5.286578186,-9.04610132)(-4.821591980,-9.20395194)
\qbezier(-4.821591980,-9.20395194)(-4.356600236,-9.36178629)(-3.882285676,-9.48888739)
\qbezier(-3.882285676,-9.48888739)(-3.407966654,-9.61597189)(-2.926354830,-9.71177921)
\qbezier(-2.926354830,-9.71177921)(-2.444739636,-9.80756966)(-1.957892882,-9.87167292)
\qbezier(-1.957892882,-9.87167292)(-1.471043868,-9.93575914)(-0.9810469362,-9.96788385)
\qbezier(-0.9810469362,-9.96788385)(-0.4910488796,-9.99999140)(0.,-10.)
\qbezier(0.,10.)(-0.4910488796,9.99999140)(-0.9810469362,9.96788385)
\qbezier(-0.9810469362,9.96788385)(-1.471043868,9.93575914)(-1.957892882,9.87167292)
\qbezier(-1.957892882,9.87167292)(-2.444739636,9.80756966)(-2.926354830,9.71177921)
\qbezier(-2.926354830,9.71177921)(-3.407966654,9.61597189)(-3.882285676,9.48888739)
\qbezier(-3.882285676,9.48888739)(-4.356600236,9.36178629)(-4.821591980,9.20395194)
\qbezier(-4.821591980,9.20395194)(-5.286578186,9.04610132)(-5.740251488,8.85819299)
\qbezier(-5.740251488,8.85819299)(-6.193918197,8.67026878)(-6.634330356,8.45309112)
\qbezier(-6.634330356,8.45309112)(-7.074734898,8.23589805)(-7.500000000,7.99038106)
\qbezier(-7.500000000,7.99038106)(-7.925256498,7.74484918)(-8.333553494,7.47204419)
\qbezier(-8.333553494,7.47204419)(-8.741840937,7.19922490)(-9.131421432,6.90030011)
\qbezier(-9.131421432,6.90030011)(-9.520991466,6.60136167)(-9.890187226,6.27759711)
\qbezier(-9.890187226,6.27759711)(-10.25937165,5.95381963)(-10.60660172,5.60660172)
\qbezier(-10.60660172,5.60660172)(-10.95381963,5.25937165)(-11.27759711,4.890187228)
\qbezier(-11.27759711,4.890187228)(-11.60136167,4.520991469)(-11.90030010,4.131421435)
\qbezier(-11.90030010,4.131421435)(-12.19922490,3.741840940)(-12.47204418,3.333553496)
\qbezier(-12.47204418,3.333553496)(-12.74484918,2.925256501)(-12.99038106,2.500000000)
\qbezier(-12.99038106,2.500000000)(-13.23589805,2.074734904)(-13.45309112,1.634330353)
\qbezier(-13.45309112,1.634330353)(-13.67026878,1.193918203)(-13.85819299,0.740251488)
\qbezier(-13.85819299,0.740251488)(-14.04610132,0.286578192)(-14.20395194,-0.178408019)
\qbezier(-14.20395194,-0.178408019)(-14.36178629,-0.643399757)(-14.48888739,-1.117714324)
\qbezier(-14.48888739,-1.117714324)(-14.61597189,-1.592033338)(-14.71177921,-2.073645168)
\qbezier(-14.71177921,-2.073645168)(-14.80756966,-2.555260355)(-14.87167292,-3.042107117)
\qbezier(-14.87167292,-3.042107117)(-14.93575914,-3.528956122)(-14.96788385,-4.018953062)
\qbezier(-14.96788385,-4.018953062)(-14.99999140,-4.508951125)(-15.,-5.)
\qbezier(-15.,-5.)(-14.99999140,-5.491048875)(-14.96788385,-5.981046938)
\qbezier(-14.96788385,-5.981046938)(-14.93575914,-6.471043878)(-14.87167292,-6.957892883)
\qbezier(-14.87167292,-6.957892883)(-14.80756966,-7.444739645)(-14.71177921,-7.926354832)
\qbezier(-14.71177921,-7.926354832)(-14.61597189,-8.407966662)(-14.48888739,-8.882285676)
\qbezier(-14.48888739,-8.882285676)(-14.36178629,-9.356600243)(-14.20395194,-9.821591981)
\qbezier(-14.20395194,-9.821591981)(-14.04610132,-10.28657819)(-13.85819299,-10.74025149)
\qbezier(-13.85819299,-10.74025149)(-13.67026878,-11.19391820)(-13.45309112,-11.63433035)
\qbezier(-13.45309112,-11.63433035)(-13.23589805,-12.07473490)(-12.99038106,-12.50000000)
\qbezier(-12.99038106,-12.50000000)(-12.74484918,-12.92525650)(-12.47204418,-13.33355350)
\qbezier(-12.47204418,-13.33355350)(-12.19922490,-13.74184094)(-11.90030010,-14.13142144)
\qbezier(-11.90030010,-14.13142144)(-11.60136167,-14.52099147)(-11.27759711,-14.89018723)
\qbezier(-11.27759711,-14.89018723)(-10.95381963,-15.25937165)(-10.60660172,-15.60660172)
\qbezier(-10.60660172,-15.60660172)(-10.25937165,-15.95381963)(-9.890187226,-16.27759711)
\qbezier(-9.890187226,-16.27759711)(-9.520991466,-16.60136167)(-9.131421432,-16.90030011)
\qbezier(-9.131421432,-16.90030011)(-8.741840937,-17.19922490)(-8.333553494,-17.47204419)
\qbezier(-8.333553494,-17.47204419)(-7.925256498,-17.74484918)(-7.500000000,-17.99038106)
\qbezier(-7.500000000,-17.99038106)(-7.074734898,-18.23589805)(-6.634330356,-18.45309112)
\qbezier(-6.634330356,-18.45309112)(-6.193918197,-18.67026878)(-5.740251488,-18.85819299)
\qbezier(-5.740251488,-18.85819299)(-5.286578186,-19.04610132)(-4.821591980,-19.20395194)
\qbezier(-4.821591980,-19.20395194)(-4.356600236,-19.36178629)(-3.882285676,-19.48888739)
\qbezier(-3.882285676,-19.48888739)(-3.407966654,-19.61597189)(-2.926354830,-19.71177921)
\qbezier(-2.926354830,-19.71177921)(-2.444739636,-19.80756966)(-1.957892882,-19.87167292)
\qbezier(-1.957892882,-19.87167292)(-1.471043868,-19.93575914)(-0.9810469362,-19.96788385)
\qbezier(-0.9810469362,-19.96788385)(-0.4910488796,-19.99999140)(0.,-20.)
\mythicklines
\put(0,-25){\line(0,1){50}}
\end{picture}
\caption{An open-book presentation of the right-handed trefoil knot}\label{fig:open-book}
\end{figure}

\medskip

Given a link $L$, let $c(L)$ denote the minimal crossing number of $L$.

\begin{theorem}[Jin-Park]\label{thm:arc nonalt}
A prime link $L$ is nonalternating if and only if
$$\alpha(L) \leq c(L).$$
\end{theorem}

\section{Kauffman polynomial}
The {\em Kauffman polynomial\/} $F_L(a,z)$ of  an oriented knot or link $L$ is defined by
$$F_L(a,z)=a^{-w(D)}\Lambda_D(a,z)$$
where $D$ is a diagram of $L$, $w(D)$ the writhe of $D$ and $\Lambda_D(a,z)$ the polynomial determined by the rules K1, K2 and K3.

\begin{itemize}
\item[(K1)] $\Lambda_O(a,z)=1$ where $O$ is the trivial knot diagram.
\item[(K2)] For any four diagrams  $D_+$, $D_-$, $D_0$ and $D_\infty$ which are identical outside a small disk in which they differ as shown below,

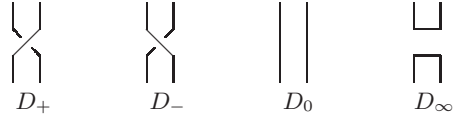
\centerline{
\small
\begin{picture}(160,50)(0,-15)
\put( 0,10){\line(0,-1){10}}\put( 0,10){\line(1,1){10}}\put(10,20){\line(0,1){10}}
\put( 0,20){\line(0,1){10}}
\qbezier( 0,20)( 0,20)( 3,17)\qbezier( 7,13)(10,10)(10,10)
\put(10,10){\line(0,-1){10}}
\put( 1,-10){$D_+$}
\put(50,10){\line(0,-1){10}}
\qbezier(50,10)(50,10)(53,13)\qbezier(57,17)(60,20)(60,20)
\put(60,20){\line(0,1){10}}
\put(50,20){\line(0,1){10}}\put(50,20){\line(1,-1){10}}\put(60,10){\line(0,-1){10}}
\put(51,-10){$D_-$}
\put(100, 0){\line(0,1){30}}\put(110, 0){\line(0,1){30}}
\put(101,-10){$D_0$}
\put(150, 0){\line(0,1){10}}\put(160, 0){\line(0,1){10}}
\put(150,10){\line(1,0){10}}\put(150,20){\line(1,0){10}}
\put(150,20){\line(0,1){10}}\put(160,20){\line(0,1){10}}
\put(150,-10){$D_\infty$}
\end{picture}
}
we have the relation
$$\Lambda_{D_+}(a,z)+\Lambda_{D_-}(a,z)=z(\Lambda_{D_0}(a,z)+\Lambda_{D_\infty}(a,z)).$$

\item[(K3)] For any three diagrams  $D_+$, $D$ and $D_-$ which are identical outside a small disk in which they differ as shown below,

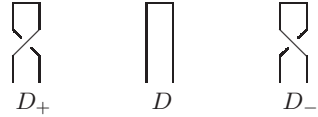
\centerline{
\small
\begin{picture}(110,50)(0,-15)
\put(0,30){\line(0,-1){10}}
\qbezier(0,20)(0,20)(3,17)\qbezier(7,13)(10,10)(10,10)
\put(10,10){\line(0,-1){10}}
\put(0,30){\line(1,0){10}}
\put(0,0){\line(0,1){10}}
\put(0,10){\line(1,1){10}}
\put(10,20){\line(0,1){10}}
\put(1,-10){$D_+$}
\put(50,30){\line(0,-1){30}}
\put(50,30){\line(1,0){10}}
\put(60,0){\line(0,1){30}}
\put(52,-10){$D$}
\put(100,30){\line(0,-1){10}}
\put(100,20){\line(1,-1){10}}
\put(110,10){\line(0,-1){10}}
\put(100,30){\line(1,0){10}}
\put(100,0){\line(0,1){10}}
\qbezier(100,10)(100,10)(103,13)\qbezier(107,17)(110,20)(110,20)
\put(110,20){\line(0,1){10}}
\put(101,-10){$D_-$}
\end{picture}
}
we have the relation
$$a\,\Lambda_{D_+}(a,z)=\Lambda_D(a,z)=a^{-1}\Lambda_{D_-}(a,z).$$
\end{itemize}

\medskip
For a connected sum and a split union of two diagrams, $\Lambda$ satisfies the following properties:
\begin{itemize}
\item[(K4)] If $D$ is a connected sum of $D_1$ and $D_2$, then
$$\Lambda_D(a,z) = \Lambda_{D_1}(a,z)\, \Lambda_{D_2}(a,z).$$
\item[(K5)] If $D$ is the split union of $D_1$ and $D_2$, then
$$\Lambda_D(a,z) = (z^{-1}a-1+z^{-1}a^{-1})\, \Lambda_{D_1}(a,z)\, \Lambda_{D_2}(a,z). $$
\end{itemize}

\medskip

The Laurent degree in the variable $a$ of the Kauffman polynomial $F_L(a,z)$ is denoted by $\operatorname{spread}_a(F_L)$ and defined by the formula				
$$\operatorname{spread}_a(F_L)= \operatorname{max-deg}_a (F_L) - \operatorname{min-deg}_a (F_L).$$
Notice that $\operatorname{spread}_a (F_L)=\operatorname{spread}_a (\Lambda_D)$ for any diagram $D$ of $L$. The following theorem gives an important lower bound for the arc index.

\begin{theorem}[Morton-Beltrami]\label{thm:lower bound of arc index}
Let $L$ be a link. Then
$$\alpha(L) \geq \operatorname{spread}_a (F_L)+2.$$
\end{theorem}

If $L$ is nonsplit and alternating, then the equality holds so that $\alpha(L)=c(L)+2$.  This is shown by Bae and Park~\cite{BP2000} using arc presentations in the form of wheel diagrams.

\section{Pretzel knots}
\begin{figure}[ht]
\centering
\setlength{\unitlength}{1.0pt}\small
\mythinlines
\begin{picture}(140,80)
\put(5,80){\line(1,0){130}}
\qbezier(5,80)(0,80)(0,75)\qbezier(5,70)(0,70)(0,75)
\qbezier(135,80)(140,80)(140,75)\qbezier(135,70)(140,70)(140,75)
\put(15,70){\line(1,0){20}}\put(45,70){\line(1,0){20}}\put(75,70){\line(1,0){15}}\put(95,69.4){$\ldots$}\put(110,70){\line(1,0){15}}
\put( 5,10){\line(0,1){5}} \put(15,10){\line(0,1){5}}
\put( 5,15){\line(1,1){10}} \qbezier(15,15)(13,17)(12,18) \qbezier( 8,22)( 7,23)( 5,25)
\put( 5,25){\line(1,1){10}} \qbezier(15,25)(13,27)(12,28) \qbezier( 8,32)( 7,33)( 5,35)
\put( 9,35){$\vdots$}\put(17,39){$p_1$}
\put( 5,45){\line(1,1){10}} \qbezier(15,45)(13,47)(12,48) \qbezier( 8,52)( 7,53)( 5,55)
\put( 5,55){\line(1,1){10}} \qbezier(15,55)(13,57)(12,58) \qbezier( 8,62)( 7,63)( 5,65)
\put( 5,65){\line(0,1){5}} \put(15,65){\line(0,1){5}}
\put(35,10){\line(0,1){5}} \put(45,10){\line(0,1){5}}
\put(35,15){\line(1,1){10}} \qbezier(45,15)(43,17)(42,18) \qbezier(38,22)(37,23)(35,25)
\put(35,25){\line(1,1){10}} \qbezier(45,25)(43,27)(42,28) \qbezier(38,32)(37,33)(35,35)
\put(39,35){$\vdots$}\put(47,39){$p_2$}
\put(35,45){\line(1,1){10}} \qbezier(45,45)(43,47)(42,48) \qbezier(38,52)(37,53)(35,55)
\put(35,55){\line(1,1){10}} \qbezier(45,55)(43,57)(42,58) \qbezier(38,62)(37,63)(35,65)
\put(35,65){\line(0,1){5}} \put(45,65){\line(0,1){5}}
\put(65,10){\line(0,1){5}} \put(75,10){\line(0,1){5}}
\put(65,65){\line(0,1){5}} \put(75,65){\line(0,1){5}}
\put(69,35){$\vdots$}
\put(95,37){$\cdots$}
\put(125,10){\line(0,1){5}} \put(135,10){\line(0,1){5}}
\put(125,15){\line(1,1){10}} \qbezier(135,15)(133,17)(132,18) \qbezier(128,22)(127,23)(125,25)
\put(125,25){\line(1,1){10}} \qbezier(135,25)(133,27)(132,28) \qbezier(128,32)(127,33)(125,35)
\put(129,35){$\vdots$}\put(137,39){$p_n$}
\put(125,45){\line(1,1){10}} \qbezier(135,45)(133,47)(132,48) \qbezier(128,52)(127,53)(125,55)
\put(125,55){\line(1,1){10}} \qbezier(135,55)(133,57)(132,58) \qbezier(128,62)(127,63)(125,65)
\put(125,65){\line(0,1){5}} \put(135,65){\line(0,1){5}}
\put(5,0){\line(1,0){130}}
\qbezier(5,0)(0,0)(0,5)\qbezier(5,10)(0,10)(0,5)
\qbezier(135,0)(140,0)(140,5)\qbezier(135,10)(140,10)(140,5)
\put(15,10){\line(1,0){20}}\put(45,10){\line(1,0){20}}\put(75,10){\line(1,0){15}}\put(95,9.5){$\ldots$}\put(110,10){\line(1,0){15}}
\end{picture}
\qquad\qquad
\mythinlines
\begin{picture}(80,80)
\put(5,80){\line(1,0){70}}
\qbezier(5,80)(0,80)(0,75)\qbezier(5,70)(0,70)(0,75)
\qbezier(75,80)(80,80)(80,75)\qbezier(75,70)(80,70)(80,75)
\put(15,70){\line(1,0){20}}\put(45,70){\line(1,0){20}}
\put( 5,10){\line(0,1){5}} \put(15,10){\line(0,1){5}}
\put(15,15){\line(-1,1){10}} \qbezier(5,15)(7,17)(8,18) \qbezier(12,22)(13,23)(15,25)
\put(15,25){\line(-1,1){10}} \qbezier(5,25)(7,27)(8,28) \qbezier(12,32)(13,33)(15,35)
\put( 9,35){$\vdots$}\put(12,39){$-p$}
\put(15,45){\line(-1,1){10}} \qbezier(5,45)(7,47)(8,48) \qbezier(12,52)(13,53)(15,55)
\put(15,55){\line(-1,1){10}} \qbezier(5,55)(7,57)(8,58) \qbezier(12,62)(13,63)(15,65)
\put( 5,65){\line(0,1){5}} \put(15,65){\line(0,1){5}}
\put(35,10){\line(0,1){5}} \put(45,10){\line(0,1){5}}
\put(35,15){\line(1,1){10}} \qbezier(45,15)(43,17)(42,18) \qbezier(38,22)(37,23)(35,25)
\put(35,25){\line(1,1){10}} \qbezier(45,25)(43,27)(42,28) \qbezier(38,32)(37,33)(35,35)
\put(39,35){$\vdots$}\put(47,39){$q$}
\put(35,45){\line(1,1){10}} \qbezier(45,45)(43,47)(42,48) \qbezier(38,52)(37,53)(35,55)
\put(35,55){\line(1,1){10}} \qbezier(45,55)(43,57)(42,58) \qbezier(38,62)(37,63)(35,65)
\put(35,65){\line(0,1){5}} \put(45,65){\line(0,1){5}}
\put(65,10){\line(0,1){5}} \put(75,10){\line(0,1){5}}
\put(65,15){\line(1,1){10}} \qbezier(75,15)(73,17)(72,18) \qbezier(68,22)(67,23)(65,25)
\put(65,25){\line(1,1){10}} \qbezier(75,25)(73,27)(72,28) \qbezier(68,32)(67,33)(65,35)
\put(69,35){$\vdots$}\put(77,39){$r$}
\put(65,45){\line(1,1){10}} \qbezier(75,45)(73,47)(72,48) \qbezier(68,52)(67,53)(65,55)
\put(65,55){\line(1,1){10}} \qbezier(75,55)(73,57)(72,58) \qbezier(68,62)(67,63)(65,65)
\put(65,65){\line(0,1){5}} \put(75,65){\line(0,1){5}}
\put(5,0){\line(1,0){70}}
\qbezier(5,0)(0,0)(0,5)\qbezier(5,10)(0,10)(0,5)
\qbezier(75,0)(80,0)(80,5)\qbezier(75,10)(80,10)(80,5)
\put(15,10){\line(1,0){20}}\put(45,10){\line(1,0){20}}
\end{picture}
\caption{Pretzel links $P(p_1,p_2,\ldots,p_n)$ and $P(-p,q,r)$}\label{fig:pretzel}
\end{figure}
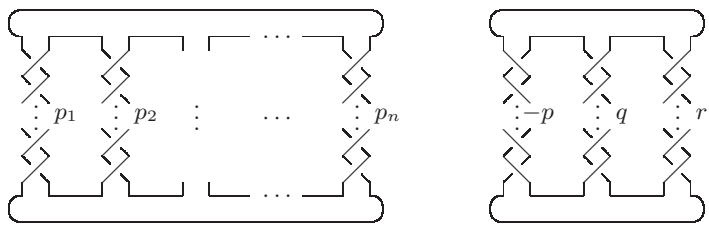

Given a sequence of integers $p_1,p_2,\ldots,p_n$, we connect two disjoint disks by $n$ bands with $p_i$ half twists, $i=1,2,\ldots,n$, so that the boundary of the resulting surface is a link as shown in Figure~\ref{fig:pretzel}. This link is called the \emph{pretzel link\/} of type $(p_1,p_2,\ldots,p_n)$ and denoted by $P(p_1,p_2,\ldots,p_n)$.

In the case $n=3$, the pretzel links satisfy the following properties:
\begin{proposition}\label{prop:pretzel_equi}
Let $p$, $q$, and $r$ be nonzero integers.
\begin{enumerate}
\item The link type of $P(p,q,r)$ is independent of the order of $p,q,r$.
\item $P(p,q,r)$ is a knot if and only if at most one of $p,q,r$ is an even number.
\end{enumerate}
\end{proposition}

\medskip

In this work, we compute the arc index for the pretzel knots $K=P(-p, q, r)$ with $p, q,r\ge2$. By Proposition~\ref{prop:pretzel_equi}(1), we may assume that $r\geq q$. By Theorem~\ref{thm:crossing}, we know that $P(-p,q,r)$ is a minimal crossing diagram of $K$, i.e., $c(K)=p+q+r$.

\begin{theorem}[Lickorish-Thistlethwaite]\label{thm:crossing}
If a link $L$ admits a reduced Montesinos diagram having $n$ crossings, then $L$ cannot be projected with fewer than $n$ crossings.
\end{theorem}

This work was motivated by Theorem~\ref{thm:upper bound of pretzel} which is a special case of Theorem~\ref{thm:arc nonalt}.

\begin{theorem}[Beltrami-Cromwell]\label{thm:upper bound of pretzel}
If $K=P(-p,q,r)$ is a knot with $p, q, r\ge2$, then
$$\alpha(K)\leq c(K)=p+q+r.$$
\end{theorem}

By computing $\operatorname{spread}_a (F_K)$ and finding arc presentations of $K=P(-p,q,r)$ with the minimum number of arcs for various values of $p$, $q$ and $r$, we obtained sharper results.

\section{Main results }

\begin{theorem}\label{thm:-2qr}
If $K=P(-2, q, r)$ is a knot with $3\le q\le r$, then
$$\alpha(K)\le c(K)-1.$$
\end{theorem}

\begin{theorem}\label{thm:-p2r}
If $K=P(-p, 2, r)$ is a knot with $p\ge 3$, $r\ge 3$, then
$$\alpha(K)=c(K).$$
\end{theorem}

\begin{theorem}\label{thm:-p3r}
If $K=P(-p, 3, r)$ is a knot with $p\ge 3$, $r\ge 3$, then
$$\alpha(K)=c(K)-1.$$
\end{theorem}

\begin{theorem}\label{thm:-p4r}
If $K=P(-p, 4, r)$ is a knot with $p\ge 5$, $r\ge 5$, then
$$\alpha(K)=c(K)-2.$$
\end{theorem}

\begin{theorem}\label{thm:-34r}
If $K=P(-3, 4, r)$ is a knot with $r\ge 7$, then
$$c(K)-4\le\alpha(K)\le c(K)-2.$$
\end{theorem}

\section{Arc presentations of $P(-p,q,r)$}

\begin{proposition}\label{prop:-2qr}
If $K=P(-2, q, r)$ is a knot with $3\le q\le r$, then $K$ has an arc presentation with $q+r+1$ arcs.
\end{proposition}

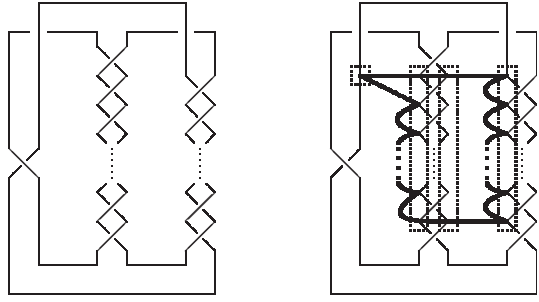
\begin{figure}[ht]
\setlength{\unitlength}{1.1pt}
\mythinlines
\begin{picture}(70,110)(0,-10)
\put(0,90){\line(1,0){7}}  \put(13,90){\line(1,0){17}} \put(40,90){\line(1,0){17}} \put(63,90){\line(1,0){7}}
\put(10,100){\line(1,0){50}}
\put(0, 0){\line(0,1){40}} \put(10,10){\line(0,1){30}}
\put(10,40){\line(-1,1){10}} \qbezier(0,40)(2,42)(4,44) \qbezier(6,46)(8,48)(10,50)
\put(0,50){\line(0,1){40}} \put(10,50){\line(0,1){50}}
\put(30,10){\line(0,1){5}} \put(40,10){\line(0,1){5}}
\put(30,15){\line(1,1){10}} \qbezier(40,15)(38,17)(36,19) \qbezier(34,21)(32,23)(30,25)
\put(30,25){\line(1,1){10}} \qbezier(40,25)(38,27)(36,29) \qbezier(34,31)(32,33)(30,35)
\qbezier(30,35)(31,36)(33,38) \qbezier(40,35)(39,36)(37,38)
\qbezier[5](35,40)(35,45)(35,50)
\qbezier(30,55)(31,54)(33,52) \qbezier(40,55)(39,54)(37,52)
\put(30,55){\line(1,1){10}} \qbezier(40,55)(38,57)(36,59) \qbezier(34,61)(32,63)(30,65)
\put(30,65){\line(1,1){10}} \qbezier(40,65)(38,67)(36,69) \qbezier(34,71)(32,73)(30,75)
\put(30,75){\line(1,1){10}} \qbezier(40,75)(38,77)(36,79) \qbezier(34,81)(32,83)(30,85)
\put(30,85){\line(0,1){5}} \put(40,85){\line(0,1){5}}
\put(60,10){\line(0,1){5}} \put(70, 0){\line(0,1){15}}
\put(60,15){\line(1,1){10}} \qbezier(70,15)(68,17)(66,19) \qbezier(64,21)(62,23)(60,25)
\put(60,25){\line(1,1){10}} \qbezier(70,25)(68,27)(66,29) \qbezier(64,31)(62,33)(60,35)
\qbezier(60,35)(61,36)(63,38) \qbezier(70,35)(69,36)(67,38)
\qbezier[5](65,40)(65,45)(65,50)
\qbezier(60,55)(61,54)(63,52) \qbezier(70,55)(69,54)(67,52)
\put(60,55){\line(1,1){10}} \qbezier(70,55)(68,57)(66,59) \qbezier(64,61)(62,63)(60,65)
\put(60,65){\line(1,1){10}} \qbezier(70,65)(68,67)(66,69) \qbezier(64,71)(62,73)(60,75)
\put(60,75){\line(0,1){25}} \put(70,75){\line(0,1){15}}
\put(10,10){\line(1,0){20}} \put(40,10){\line(1,0){20}}
\put(0,0){\line(1,0){70}}
\end{picture}
\qquad\qquad
\begin{picture}(70,110)(0,-10)
\put(0,90){\line(1,0){7}}  \put(13,90){\line(1,0){17}} \put(40,90){\line(1,0){17}} \put(63,90){\line(1,0){7}}
\put(10,100){\line(1,0){50}}
\put(0, 0){\line(0,1){40}} \put(10,10){\line(0,1){30}}
\put(10,40){\line(-1,1){10}} \qbezier(0,40)(2,42)(4,44) \qbezier(6,46)(8,48)(10,50)
\put(0,50){\line(0,1){40}} \put(10,50){\line(0,1){50}}
\put(30,10){\line(0,1){5}} \put(40,10){\line(0,1){5}}
\put(30,15){\line(1,1){10}} \qbezier(40,15)(38,17)(36,19) \qbezier(34,21)(32,23)(30,25)
\put(30,25){\line(1,1){10}} \qbezier(40,25)(38,27)(36,29) \qbezier(34,31)(32,33)(30,35)
\qbezier(30,35)(31,36)(33,38) \qbezier(40,35)(39,36)(37,38)
\qbezier[5](35,40)(35,45)(35,50)
\qbezier(30,55)(31,54)(33,52) \qbezier(40,55)(39,54)(37,52)
\put(30,55){\line(1,1){10}} \qbezier(40,55)(38,57)(36,59) \qbezier(34,61)(32,63)(30,65)
\put(30,65){\line(1,1){10}} \qbezier(40,65)(38,67)(36,69) \qbezier(34,71)(32,73)(30,75)
\put(30,75){\line(1,1){10}} \qbezier(40,75)(38,77)(36,79) \qbezier(34,81)(32,83)(30,85)
\put(30,85){\line(0,1){5}} \put(40,85){\line(0,1){5}}
\put(60,10){\line(0,1){5}} \put(70, 0){\line(0,1){15}}
\put(60,15){\line(1,1){10}} \qbezier(70,15)(68,17)(66,19) \qbezier(64,21)(62,23)(60,25)
\put(60,25){\line(1,1){10}} \qbezier(70,25)(68,27)(66,29) \qbezier(64,31)(62,33)(60,35)
\qbezier(60,35)(61,36)(63,38) \qbezier(70,35)(69,36)(67,38)
\qbezier[5](65,40)(65,45)(65,50)
\qbezier(60,55)(61,54)(63,52) \qbezier(70,55)(69,54)(67,52)
\put(60,55){\line(1,1){10}} \qbezier(70,55)(68,57)(66,59) \qbezier(64,61)(62,63)(60,65)
\put(60,65){\line(1,1){10}} \qbezier(70,65)(68,67)(66,69) \qbezier(64,71)(62,73)(60,75)
\put(60,75){\line(0,1){25}} \put(70,75){\line(0,1){15}}
\put(10,10){\line(1,0){20}} \put(40,10){\line(1,0){20}}
\put(0,0){\line(1,0){70}}
\mythicklines
\qbezier(10,75)(35,75)(60,75)
\qbezier(60,75)(45,70)(60,65) \qbezier(60,65)(45,60)(60,55) \qbezier(60,55)(55,55)(53,52)
\qbezier[3](53,40)(53,45)(53,50)
\qbezier(60,35)(55,35)(53,38) \qbezier(60,25)(45,30)(60,35)
\qbezier(35,25)(15,25)(30,35) \qbezier(35,25)(45,25)(60,25)
\qbezier(30,35)(25,35)(23,38)
\qbezier(10,75)(20,70)(30,65) \qbezier(30,65)(15,60)(30,55) \qbezier(30,55)(25,55)(23,52)
\qbezier[3](23,40)(23,45)(23,50)
\thicklines
\qbezier[4](7,78)(10,78)(13,78)\qbezier[4](7,72)(10,72)(13,72)\qbezier[4](7,72)(7,75)(7,78)\qbezier[4](13,72)(13,75)(13,78)
\qbezier[4](27,78)(30,78)(33,78)\qbezier[4](27,22)(30,22)(33,22)\qbezier[38](27,22)(27,50)(27,78)\qbezier[38](33,22)(33,50)(33,78)
\qbezier[4](37,78)(40,78)(43,78)\qbezier[4](37,22)(40,22)(43,22)\qbezier[38](37,22)(37,50)(37,78)\qbezier[38](43,22)(43,50)(43,78)
\qbezier[4](57,78)(60,78)(63,78)\qbezier[4](57,22)(60,22)(63,22)\qbezier[38](57,22)(57,50)(57,78)\qbezier[38](63,22)(63,50)(63,78)
\end{picture}
\vspace{-4mm}
\caption{An arc presentation of $P(-2,q,r)$}\label{fig:(-2,q,r)}
\end{figure}

\begin{proof}
Figure~\ref{fig:(-2,q,r)} shows a pretzel diagram of  $P(-2,q,r)$ and its arc presentation with $q+r+1$ arcs. The thick curve is the axis of the arc presentation which cuts the knot at 1 place in the leftmost box, $q-1$ places in the second, 2 places in the third, and $r-1$ places in the fourth. The $q+r+1$ arcs of the knot satisfies the four properties of an arc presentation.
\end{proof}

\begin{proposition}\label{prop:-pqr,q=2,3}
If $K=P(-p,q,r)$ is a knot with $p\ge3$ and $2\le q\le 3\le r$, then
then $K$ has an arc presentation with $p+r+2$ arcs.
\end{proposition}

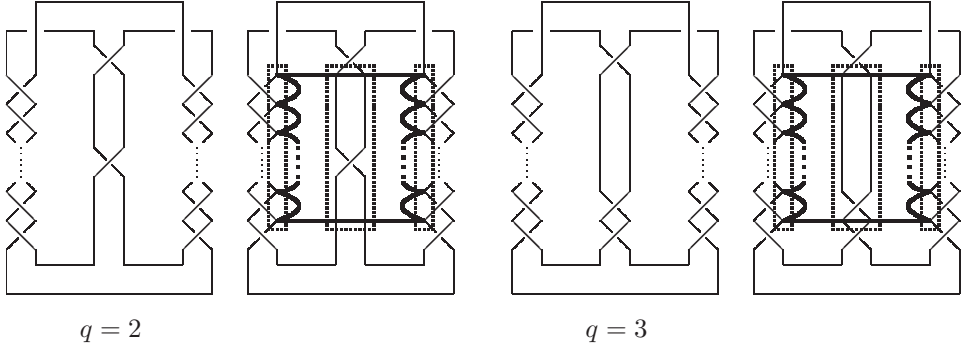
\begin{figure}[ht]
\setlength{\unitlength}{1.1pt}
\mythinlines
\begin{picture}(70,120)(0,-20)
\put(0,90){\line(1,0){7}} \put(13,90){\line(1,0){17}} \put(40,90){\line(1,0){17}} \put(63,90){\line(1,0){7}}
\put(10,100){\line(1,0){50}}
\put(0, 0){\line(0,1){15}} \put(10,10){\line(0,1){5}}
\put(10,15){\line(-1,1){10}} \qbezier(0,15)(2,17)(4,19) \qbezier(6,21)(8,23)(10,25)
\put(10,25){\line(-1,1){10}} \qbezier(0,25)(2,27)(4,29) \qbezier(6,31)(8,33)(10,35)
\put(0,75){\line(0,1){15}} \put(10,75){\line(0,1){25}}
\qbezier(10,35)(9,36)(7,38) \qbezier(0,35)(1,36)(3,38)
\qbezier[5](5,40)(5,45)(5,50)
\qbezier(10,55)(9,54)(7,52) \qbezier(0,55)(1,54)(3,52)
\put(10,55){\line(-1,1){10}} \qbezier(0,55)(2,57)(4,59) \qbezier(6,61)(8,63)(10,65)
\put(10,65){\line(-1,1){10}} \qbezier(0,65)(2,67)(4,69) \qbezier(6,71)(8,73)(10,75)
\put(30,10){\line(0,1){30}} \put(40,10){\line(0,1){30}}
\put(30,40){\line(1,1){10}} \qbezier(40,40)(38,42)(36,44) \qbezier(34,46)(32,48)(30,50)
\put(30,50){\line(0,1){25}} \put(40,50){\line(0,1){25}}
\put(30,75){\line(1,1){10}} \qbezier(40,75)(38,77)(36,79) \qbezier(34,81)(32,83)(30,85)
\put(30,85){\line(0,1){5}} \put(40,85){\line(0,1){5}}
\put(60,10){\line(0,1){5}} \put(70, 0){\line(0,1){15}}
\put(60,15){\line(1,1){10}} \qbezier(70,15)(68,17)(66,19) \qbezier(64,21)(62,23)(60,25)
\put(60,25){\line(1,1){10}} \qbezier(70,25)(68,27)(66,29) \qbezier(64,31)(62,33)(60,35)
\qbezier(60,35)(61,36)(63,38) \qbezier(70,35)(69,36)(67,38)
\qbezier[5](65,40)(65,45)(65,50)
\qbezier(60,55)(61,54)(63,52) \qbezier(70,55)(69,54)(67,52)
\put(60,55){\line(1,1){10}} \qbezier(70,55)(68,57)(66,59) \qbezier(64,61)(62,63)(60,65)
\put(60,65){\line(1,1){10}} \qbezier(70,65)(68,67)(66,69) \qbezier(64,71)(62,73)(60,75)
\put(60,75){\line(0,1){25}} \put(70,75){\line(0,1){15}}
\put(10,10){\line(1,0){20}} \put(40,10){\line(1,0){20}}
\put(0,0){\line(1,0){70}}
\put(0,-15){\hbox to 78pt{\hfill$q=2$\hfill}}
\end{picture}
\quad
\begin{picture}(70,120)(0,-20)
\put(0,90){\line(1,0){7}} \put(13,90){\line(1,0){17}} \put(40,90){\line(1,0){17}} \put(63,90){\line(1,0){7}}
\put(10,100){\line(1,0){50}}
\put(0, 0){\line(0,1){15}} \put(10,10){\line(0,1){5}}
\put(10,15){\line(-1,1){10}} \qbezier(0,15)(2,17)(4,19) \qbezier(6,21)(8,23)(10,25)
\put(10,25){\line(-1,1){10}} \qbezier(0,25)(2,27)(4,29) \qbezier(6,31)(8,33)(10,35)
\put(0,75){\line(0,1){15}} \put(10,75){\line(0,1){25}}
\qbezier(10,35)(9,36)(7,38) \qbezier(0,35)(1,36)(3,38)
\qbezier[5](5,40)(5,45)(5,50)
\qbezier(10,55)(9,54)(7,52) \qbezier(0,55)(1,54)(3,52)
\put(10,55){\line(-1,1){10}} \qbezier(0,55)(2,57)(4,59) \qbezier(6,61)(8,63)(10,65)
\put(10,65){\line(-1,1){10}} \qbezier(0,65)(2,67)(4,69) \qbezier(6,71)(8,73)(10,75)
\put(30,10){\line(0,1){30}} \put(40,10){\line(0,1){30}}
\put(30,40){\line(1,1){10}} \qbezier(40,40)(38,42)(36,44) \qbezier(34,46)(32,48)(30,50)
\put(30,50){\line(0,1){25}} \put(40,50){\line(0,1){25}}
\put(30,75){\line(1,1){10}} \qbezier(40,75)(38,77)(36,79) \qbezier(34,81)(32,83)(30,85)
\put(30,85){\line(0,1){5}} \put(40,85){\line(0,1){5}}
\put(60,10){\line(0,1){5}} \put(70, 0){\line(0,1){15}}
\put(60,15){\line(1,1){10}} \qbezier(70,15)(68,17)(66,19) \qbezier(64,21)(62,23)(60,25)
\put(60,25){\line(1,1){10}} \qbezier(70,25)(68,27)(66,29) \qbezier(64,31)(62,33)(60,35)
\qbezier(60,35)(61,36)(63,38) \qbezier(70,35)(69,36)(67,38)
\qbezier[5](65,40)(65,45)(65,50)
\qbezier(60,55)(61,54)(63,52) \qbezier(70,55)(69,54)(67,52)
\put(60,55){\line(1,1){10}} \qbezier(70,55)(68,57)(66,59) \qbezier(64,61)(62,63)(60,65)
\put(60,65){\line(1,1){10}} \qbezier(70,65)(68,67)(66,69) \qbezier(64,71)(62,73)(60,75)
\put(60,75){\line(0,1){25}} \put(70,75){\line(0,1){15}}
\put(10,10){\line(1,0){20}} \put(40,10){\line(1,0){20}}
\put(0,0){\line(1,0){70}}
\mythicklines
\qbezier(10,75)(35,75)(60,75)
\qbezier(60,75)(45,70)(60,65) \qbezier(60,65)(45,60)(60,55) \qbezier(60,55)(55,55)(53,52)
\qbezier[3](53,40)(53,45)(53,50)
\qbezier(60,35)(55,35)(53,38) \qbezier(60,25)(45,30)(60,35)
\put(10,25){\line(1,0){50}}
\qbezier(10,35)(15,35)(17,38) \qbezier(10,25)(25,30)(10,35)
\qbezier(10,75)(25,70)(10,65) \qbezier(10,65)(25,60)(10,55) \qbezier(10,55)(15,55)(17,52)
\qbezier[3](17,40)(17,45)(17,50)
\thicklines
\qbezier[4](7,78)(10,78)(13,78)\qbezier[4](7,22)(10,22)(13,22)\qbezier[38](7,22)(7,50)(7,78)\qbezier[38](13,22)(13,50)(13,78)
\qbezier[10](27,78)(35,78)(43,78)\qbezier[10](27,22)(35,22)(43,22)\qbezier[38](27,22)(27,50)(27,78)\qbezier[38](43,22)(43,50)(43,78)
\qbezier[4](57,78)(60,78)(63,78)\qbezier[4](57,22)(60,22)(63,22)\qbezier[38](57,22)(57,50)(57,78)\qbezier[38](63,22)(63,50)(63,78)
\end{picture}
\hfill
\mythinlines
\begin{picture}(70,120)(0,-20)
\put(0,90){\line(1,0){7}} \put(13,90){\line(1,0){17}} \put(40,90){\line(1,0){17}} \put(63,90){\line(1,0){7}}
\put(10,100){\line(1,0){50}}
\put(0, 0){\line(0,1){15}} \put(10,10){\line(0,1){5}}
\put(10,15){\line(-1,1){10}} \qbezier(0,15)(2,17)(4,19) \qbezier(6,21)(8,23)(10,25)
\put(10,25){\line(-1,1){10}} \qbezier(0,25)(2,27)(4,29) \qbezier(6,31)(8,33)(10,35)
\put(0,75){\line(0,1){15}} \put(10,75){\line(0,1){25}}
\qbezier(10,35)(9,36)(7,38) \qbezier(0,35)(1,36)(3,38)
\qbezier[5](5,40)(5,45)(5,50)
\qbezier(10,55)(9,54)(7,52) \qbezier(0,55)(1,54)(3,52)
\put(10,55){\line(-1,1){10}} \qbezier(0,55)(2,57)(4,59) \qbezier(6,61)(8,63)(10,65)
\put(10,65){\line(-1,1){10}} \qbezier(0,65)(2,67)(4,69) \qbezier(6,71)(8,73)(10,75)
\put(30,10){\line(0,1){5}} \put(40,10){\line(0,1){5}}
\put(30,15){\line(1,1){10}} \qbezier(40,15)(38,17)(36,19) \qbezier(34,21)(32,23)(30,25)
\put(30,25){\line(1,1){10}} \qbezier(40,25)(38,27)(36,29) \qbezier(34,31)(32,33)(30,35)
\put(30,35){\line(0,1){40}} \put(40,35){\line(0,1){40}}
\put(30,75){\line(1,1){10}} \qbezier(40,75)(38,77)(36,79) \qbezier(34,81)(32,83)(30,85)
\put(30,85){\line(0,1){5}} \put(40,85){\line(0,1){5}}
\put(60,10){\line(0,1){5}} \put(70, 0){\line(0,1){15}}
\put(60,15){\line(1,1){10}} \qbezier(70,15)(68,17)(66,19) \qbezier(64,21)(62,23)(60,25)
\put(60,25){\line(1,1){10}} \qbezier(70,25)(68,27)(66,29) \qbezier(64,31)(62,33)(60,35)
\qbezier(60,35)(61,36)(63,38) \qbezier(70,35)(69,36)(67,38)
\qbezier[5](65,40)(65,45)(65,50)
\qbezier(60,55)(61,54)(63,52) \qbezier(70,55)(69,54)(67,52)
\put(60,55){\line(1,1){10}} \qbezier(70,55)(68,57)(66,59) \qbezier(64,61)(62,63)(60,65)
\put(60,65){\line(1,1){10}} \qbezier(70,65)(68,67)(66,69) \qbezier(64,71)(62,73)(60,75)
\put(60,75){\line(0,1){25}} \put(70,75){\line(0,1){15}}
\put(10,10){\line(1,0){20}} \put(40,10){\line(1,0){20}}
\put(0,0){\line(1,0){70}}
\put(0,-15){\hbox to 78pt{\hfill$q=3$\hfill}}
\end{picture}
\quad
\begin{picture}(70,120)(0,-20)
\put(0,90){\line(1,0){7}} \put(13,90){\line(1,0){17}} \put(40,90){\line(1,0){17}} \put(63,90){\line(1,0){7}}
\put(10,100){\line(1,0){50}}
\put(0, 0){\line(0,1){15}} \put(10,10){\line(0,1){5}}
\put(10,15){\line(-1,1){10}} \qbezier(0,15)(2,17)(4,19) \qbezier(6,21)(8,23)(10,25)
\put(10,25){\line(-1,1){10}} \qbezier(0,25)(2,27)(4,29) \qbezier(6,31)(8,33)(10,35)
\put(0,75){\line(0,1){15}} \put(10,75){\line(0,1){25}}
\qbezier(10,35)(9,36)(7,38) \qbezier(0,35)(1,36)(3,38)
\qbezier[5](5,40)(5,45)(5,50)
\qbezier(10,55)(9,54)(7,52) \qbezier(0,55)(1,54)(3,52)
\put(10,55){\line(-1,1){10}} \qbezier(0,55)(2,57)(4,59) \qbezier(6,61)(8,63)(10,65)
\put(10,65){\line(-1,1){10}} \qbezier(0,65)(2,67)(4,69) \qbezier(6,71)(8,73)(10,75)
\put(30,10){\line(0,1){5}} \put(40,10){\line(0,1){5}}
\put(30,15){\line(1,1){10}} \qbezier(40,15)(38,17)(36,19) \qbezier(34,21)(32,23)(30,25)
\put(30,25){\line(1,1){10}} \qbezier(40,25)(38,27)(36,29) \qbezier(34,31)(32,33)(30,35)
\put(30,35){\line(0,1){40}} \put(40,35){\line(0,1){40}}
\put(30,75){\line(1,1){10}} \qbezier(40,75)(38,77)(36,79) \qbezier(34,81)(32,83)(30,85)
\put(30,85){\line(0,1){5}} \put(40,85){\line(0,1){5}}
\put(60,10){\line(0,1){5}} \put(70, 0){\line(0,1){15}}
\put(60,15){\line(1,1){10}} \qbezier(70,15)(68,17)(66,19) \qbezier(64,21)(62,23)(60,25)
\put(60,25){\line(1,1){10}} \qbezier(70,25)(68,27)(66,29) \qbezier(64,31)(62,33)(60,35)
\qbezier(60,35)(61,36)(63,38) \qbezier(70,35)(69,36)(67,38)
\qbezier[5](65,40)(65,45)(65,50)
\qbezier(60,55)(61,54)(63,52) \qbezier(70,55)(69,54)(67,52)
\put(60,55){\line(1,1){10}} \qbezier(70,55)(68,57)(66,59) \qbezier(64,61)(62,63)(60,65)
\put(60,65){\line(1,1){10}} \qbezier(70,65)(68,67)(66,69) \qbezier(64,71)(62,73)(60,75)
\put(60,75){\line(0,1){25}} \put(70,75){\line(0,1){15}}
\put(10,10){\line(1,0){20}} \put(40,10){\line(1,0){20}}
\put(0,0){\line(1,0){70}}
\mythicklines
\qbezier(10,75)(35,75)(60,75)
\qbezier(60,75)(45,70)(60,65) \qbezier(60,65)(45,60)(60,55) \qbezier(60,55)(55,55)(53,52)
\qbezier[3](53,40)(53,45)(53,50)
\qbezier(60,35)(55,35)(53,38) \qbezier(60,25)(45,30)(60,35)
\qbezier(10,25)(35,25)(60,25)
\qbezier(10,35)(15,35)(17,38) \qbezier(10,25)(25,30)(10,35)
\qbezier(10,75)(25,70)(10,65) \qbezier(10,65)(25,60)(10,55) \qbezier(10,55)(15,55)(17,52)
\qbezier[3](17,40)(17,45)(17,50)
\thicklines
\qbezier[4](7,78)(10,78)(13,78)\qbezier[4](7,22)(10,22)(13,22)\qbezier[38](7,22)(7,50)(7,78)\qbezier[38](13,22)(13,50)(13,78)
\qbezier[10](27,78)(35,78)(43,78)\qbezier[10](27,22)(35,22)(43,22)\qbezier[38](27,22)(27,50)(27,78)\qbezier[38](43,22)(43,50)(43,78)
\qbezier[4](57,78)(60,78)(63,78)\qbezier[4](57,22)(60,22)(63,22)\qbezier[38](57,22)(57,50)(57,78)\qbezier[38](63,22)(63,50)(63,78)
\end{picture}
\vspace{-4mm}
\caption{Arc presentations of $P(-p,q,r)$ for $q=2,3$}\label{fig:(-p,2 or 3,r)}
\end{figure}

\begin{proof}
For each of $q=2,3$, Figure~\ref{fig:(-p,2 or 3,r)} shows a pretzel diagram of  $P(-p,q,r)$ and its arc presentation with $p+r+2$ arcs. The thick curve is the axis of the arc presentation which cuts the knot at $p-1$ places in the leftmost box, $4$ places in the second, and $r-1$ places in the third. The $p+r+2$ arcs of the knot satisfies the four properties of an arc presentation.
\end{proof}

\begin{proposition}\label{prop:-pqr,q>3}
If $K=P(-p, q, r)$ is a knot with $p\geq3$ and  $4\le q\le r$, then $K$ has an arc presentation with $p+q+r-2$ arcs.
\end{proposition}

\begin{proof}
In Figure~\ref{fig:(-p,q,r)}, the diagram (a) shows a pretzel diagram of $P(-p,q,r)$ with $p\geq3$ and  $4\le q\le r$. The diagram (b) is obtained from (a) by two applications of the Reidemeister move of type 3. The diagram (c) shows an arc presentation with $p+q+r-1$ arcs.
The diagram (d) is obtained from (c) by isotoping the arc labeled $x$ over the axis so that there are only $p+q+r-2$ arcs. Each of the seven boxes, from left to right, contains $1$, $p-3$, $2$, $3$, $q-3$, $1$, and $r-3$ arcs, respectively.
\end{proof}

\begin{figure}[hb]
\setlength{\unitlength}{1.1pt}
\mythinlines
\begin{picture}(90,115)(-15,-15)
\put(0,90){\line(1,0){7}} \put(13,90){\line(1,0){17}} \put(40,90){\line(1,0){17}} \put(63,90){\line(1,0){7}}
\put(10,100){\line(1,0){50}}
\put(0,0){\line(0,1){15}} \put(10,-10){\line(0,1){7}} \put(10,3){\line(0,1){12}}
\put(10,15){\line(-1,1){10}} \qbezier(0,15)(2,17)(4,19) \qbezier(6,21)(8,23)(10,25)
\put(10,25){\line(-1,1){10}} \qbezier(0,25)(2,27)(4,29) \qbezier(6,31)(8,33)(10,35)
\put(0,75){\line(0,1){15}} \put(10,75){\line(0,1){25}}
\qbezier(10,35)(9,36)(7,38) \qbezier(0,35)(1,36)(3,38)
\qbezier[5](5,40)(5,45)(5,50)
\qbezier(10,55)(9,54)(7,52) \qbezier(0,55)(1,54)(3,52)
\put(10,55){\line(-1,1){10}} \qbezier(0,55)(2,57)(4,59) \qbezier(6,61)(8,63)(10,65)
\put(10,65){\line(-1,1){10}} \qbezier(0,65)(2,67)(4,69) \qbezier(6,71)(8,73)(10,75)
\put(30,0){\line(0,1){5}} \put(40,0){\line(0,1){5}}
\put(30,5){\line(1,1){10}} \qbezier(40,5)(38,7)(36,9) \qbezier(34,11)(32,13)(30,15)
\put(30,15){\line(1,1){10}} \qbezier(40,15)(38,17)(36,19) \qbezier(34,21)(32,23)(30,25)
\put(30,25){\line(1,1){10}} \qbezier(40,25)(38,27)(36,29) \qbezier(34,31)(32,33)(30,35)
\qbezier(30,35)(31,36)(33,38) \qbezier(40,35)(39,36)(37,38)
\qbezier[5](35,40)(35,45)(35,50)
\qbezier(30,55)(31,54)(33,52) \qbezier(40,55)(39,54)(37,52)
\put(30,55){\line(1,1){10}} \qbezier(40,55)(38,57)(36,59) \qbezier(34,61)(32,63)(30,65)
\put(30,65){\line(1,1){10}} \qbezier(40,65)(38,67)(36,69) \qbezier(34,71)(32,73)(30,75)
\put(30,75){\line(1,1){10}} \qbezier(40,75)(38,77)(36,79) \qbezier(34,81)(32,83)(30,85)
\put(30,85){\line(0,1){5}} \put(40,85){\line(0,1){5}}
\put(60,-10){\line(0,1){7}} \put(60,3){\line(0,1){12}} \put(70,0){\line(0,1){15}}
\put(60,15){\line(1,1){10}} \qbezier(70,15)(68,17)(66,19) \qbezier(64,21)(62,23)(60,25)
\put(60,25){\line(1,1){10}} \qbezier(70,25)(68,27)(66,29) \qbezier(64,31)(62,33)(60,35)
\qbezier(60,35)(61,36)(63,38) \qbezier(70,35)(69,36)(67,38)
\qbezier[5](65,40)(65,45)(65,50)
\qbezier(60,55)(61,54)(63,52) \qbezier(70,55)(69,54)(67,52)
\put(60,55){\line(1,1){10}} \qbezier(70,55)(68,57)(66,59) \qbezier(64,61)(62,63)(60,65)
\put(60,65){\line(1,1){10}} \qbezier(70,65)(68,67)(66,69) \qbezier(64,71)(62,73)(60,75)
\put(60,75){\line(0,1){25}} \put(70,75){\line(0,1){15}}
\put(0,0){\line(1,0){30}} \put(40,0){\line(1,0){30}}
\put(10,-10){\line(1,0){50}}
\put(-15,82){(a)}
\end{picture}
\qquad
\begin{picture}(90,115)(-15,-15)
\put(0,100){\line(1,0){30}} \put(40,100){\line(1,0){30}}
\put(10,80){\line(1,0){50}}
\put(0,-10){\line(0,1){25}} \put(10,10){\line(0,1){5}}
\put(10,15){\line(-1,1){10}} \qbezier(0,15)(2,17)(4,19) \qbezier(6,21)(8,23)(10,25)
\put(10,25){\line(-1,1){10}} \qbezier(0,25)(2,27)(4,29) \qbezier(6,31)(8,33)(10,35)
\put(0,75){\line(0,1){25}} \put(10,75){\line(0,1){5}}
\qbezier(10,35)(9,36)(7,38) \qbezier(0,35)(1,36)(3,38)
\qbezier[5](5,40)(5,45)(5,50)
\qbezier(10,55)(9,54)(7,52) \qbezier(0,55)(1,54)(3,52)
\put(10,55){\line(-1,1){10}} \qbezier(0,55)(2,57)(4,59) \qbezier(6,61)(8,63)(10,65)
\put(10,65){\line(-1,1){10}} \qbezier(0,65)(2,67)(4,69) \qbezier(6,71)(8,73)(10,75)
\put(30,-10){\line(0,1){5}} \put(40,-10){\line(0,1){5}}
\put(30,-5){\line(1,1){10}} \qbezier(40,-5)(38,-3)(36,-1) \qbezier(34,1)(32,3)(30,5)
\put(30,5){\line(0,1){10}} \put(40,5){\line(0,1){10}}
\put(30,15){\line(1,1){10}} \qbezier(40,15)(38,17)(36,19) \qbezier(34,21)(32,23)(30,25)
\put(30,25){\line(1,1){10}} \qbezier(40,25)(38,27)(36,29) \qbezier(34,31)(32,33)(30,35)
\qbezier(30,35)(31,36)(33,38) \qbezier(40,35)(39,36)(37,38)
\qbezier[5](35,40)(35,45)(35,50)
\qbezier(30,55)(31,54)(33,52) \qbezier(40,55)(39,54)(37,52)
\put(30,55){\line(1,1){10}} \qbezier(40,55)(38,57)(36,59) \qbezier(34,61)(32,63)(30,65)
\put(30,65){\line(1,1){10}} \qbezier(40,65)(38,67)(36,69) \qbezier(34,71)(32,73)(30,75)
\put(30,75){\line(0,1){3}} \put(40,75){\line(0,1){3}}
\put(30,82){\line(0,1){3}} \put(40,82){\line(0,1){3}}
\put(30,85){\line(1,1){10}} \qbezier(40,85)(38,87)(36,89) \qbezier(34,91)(32,93)(30,95)
\put(30,95){\line(0,1){5}} \put(40,95){\line(0,1){5}}
\put(60,10){\line(0,1){5}} \put(70,-10){\line(0,1){25}}
\put(60,15){\line(1,1){10}} \qbezier(70,15)(68,17)(66,19) \qbezier(64,21)(62,23)(60,25)
\put(60,25){\line(1,1){10}} \qbezier(70,25)(68,27)(66,29) \qbezier(64,31)(62,33)(60,35)
\qbezier(60,35)(61,36)(63,38) \qbezier(70,35)(69,36)(67,38)
\qbezier[5](65,40)(65,45)(65,50)
\qbezier(60,55)(61,54)(63,52) \qbezier(70,55)(69,54)(67,52)
\put(60,55){\line(1,1){10}} \qbezier(70,55)(68,57)(66,59) \qbezier(64,61)(62,63)(60,65)
\put(60,65){\line(1,1){10}} \qbezier(70,65)(68,67)(66,69) \qbezier(64,71)(62,73)(60,75)
\put(60,75){\line(0,1){5}} \put(70,75){\line(0,1){25}}
\put(0,-10){\line(1,0){30}} \put(40,-10){\line(1,0){30}}
\put(10,10){\line(1,0){17}} \put(33,10){\line(1,0){4}} \put(43,10){\line(1,0){17}}
\put(-15,82){(b)}
\end{picture}

\begin{picture}(90,135)(-15,-20)
\put(0,100){\line(1,0){30}} \put(40,100){\line(1,0){30}}
\put(10,80){\line(1,0){50}}
\put(0,-10){\line(0,1){25}} \put(10,10){\line(0,1){5}}
\put(10,15){\line(-1,1){10}} \qbezier(0,15)(2,17)(4,19) \qbezier(6,21)(8,23)(10,25)
\put(10,25){\line(-1,1){10}} \qbezier(0,25)(2,27)(4,29) \qbezier(6,31)(8,33)(10,35)
\put(0,75){\line(0,1){25}} \put(10,75){\line(0,1){5}}
\qbezier(10,35)(9,36)(7,38) \qbezier(0,35)(1,36)(3,38)
\qbezier[5](5,40)(5,45)(5,50)
\qbezier(10,55)(9,54)(7,52) \qbezier(0,55)(1,54)(3,52)
\put(10,55){\line(-1,1){10}} \qbezier(0,55)(2,57)(4,59) \qbezier(6,61)(8,63)(10,65)
\put(10,65){\line(-1,1){10}} \qbezier(0,65)(2,67)(4,69) \qbezier(6,71)(8,73)(10,75)
\put(30,-10){\line(0,1){5}} \put(40,-10){\line(0,1){5}}
\put(30,-5){\line(1,1){10}} \qbezier(40,-5)(38,-3)(36,-1) \qbezier(34,1)(32,3)(30,5)
\put(30,5){\line(0,1){10}} \put(40,5){\line(0,1){10}}
\put(30,15){\line(1,1){10}} \qbezier(40,15)(38,17)(36,19) \qbezier(34,21)(32,23)(30,25)
\put(30,25){\line(1,1){10}} \qbezier(40,25)(38,27)(36,29) \qbezier(34,31)(32,33)(30,35)
\qbezier(30,35)(31,36)(33,38) \qbezier(40,35)(39,36)(37,38)
\qbezier[5](35,40)(35,45)(35,50)
\qbezier(30,55)(31,54)(33,52) \qbezier(40,55)(39,54)(37,52)
\put(30,55){\line(1,1){10}} \qbezier(40,55)(38,57)(36,59) \qbezier(34,61)(32,63)(30,65)
\put(30,65){\line(1,1){10}} \qbezier(40,65)(38,67)(36,69) \qbezier(34,71)(32,73)(30,75)
\put(30,75){\line(0,1){3}} \put(40,75){\line(0,1){3}}
\put(30,82){\line(0,1){3}} \put(40,82){\line(0,1){3}}
\put(30,85){\line(1,1){10}} \qbezier(40,85)(38,87)(36,89) \qbezier(34,91)(32,93)(30,95)
\put(30,95){\line(0,1){5}} \put(40,95){\line(0,1){5}}
\put(60,10){\line(0,1){5}} \put(70,-10){\line(0,1){25}}
\put(60,15){\line(1,1){10}} \qbezier(70,15)(68,17)(66,19) \qbezier(64,21)(62,23)(60,25)
\put(60,25){\line(1,1){10}} \qbezier(70,25)(68,27)(66,29) \qbezier(64,31)(62,33)(60,35)
\qbezier(60,35)(61,36)(63,38) \qbezier(70,35)(69,36)(67,38)
\qbezier[5](65,40)(65,45)(65,50)
\qbezier(60,55)(61,54)(63,52) \qbezier(70,55)(69,54)(67,52)
\put(60,55){\line(1,1){10}} \qbezier(70,55)(68,57)(66,59) \qbezier(64,61)(62,63)(60,65)
\put(60,65){\line(1,1){10}} \qbezier(70,65)(68,67)(66,69) \qbezier(64,71)(62,73)(60,75)
\put(60,75){\line(0,1){5}} \put(70,75){\line(0,1){25}}
\put(0,-10){\line(1,0){30}} \put(40,-10){\line(1,0){30}}
\put(10,10){\line(1,0){17}} \put(33,10){\line(1,0){4}} \put(43,10){\line(1,0){17}}
\mythicklines
\qbezier(10,65)(20,65)(40,65)
\qbezier(40,65)(55,60)(40,55) \qbezier(40,55)(45,55)(47,52)
\qbezier[3](47,40)(47,45)(47,50)
\qbezier(40,35)(45,35)(47,38) \qbezier(40,25)(55,30)(40,35)
\qbezier(40,25)(50,25)(70,25)
\qbezier(70,25)(85,30)(70,35) \qbezier(70,35)(75,35)(77,38)
\qbezier[3](77,40)(77,45)(77,50)
\qbezier(70,55)(85,60)(70,65) \qbezier(70,55)(75,55)(77,52)
\qbezier(70,65)(100,115)(50,115) \qbezier(50,115)(0,115)(30,85)
\qbezier(30,85)(-10,105)(-15,45)
\qbezier(-15,45)(-10,-15)(30,5)
\qbezier(30,5)(22,5)(20,10) \qbezier(10,25)(18,25)(20,10)
\qbezier(10,35)(15,35)(17,38) \qbezier(10,25)(25,30)(10,35)
\qbezier(10,65)(25,60)(10,55) \qbezier(10,55)(15,55)(17,52)
\qbezier[3](17,40)(17,45)(17,50)
\put(2,94){$x$}
\put(-15,92){(c)}
\end{picture}
\qquad
\begin{picture}(90,135)(-15,-20)
\put(20,105){\line(1,0){30}}
\put(40,100){\line(1,0){8}} \put(52,100){\line(1,0){18}}
\put(10,79){\line(1,0){50}}
\put(0,81){\line(1,0){50}}
\put(0,-10){\line(0,1){25}} \put(10,10){\line(0,1){5}}
\put(10,15){\line(-1,1){10}} \qbezier(0,15)(2,17)(4,19) \qbezier(6,21)(8,23)(10,25)
\put(10,25){\line(-1,1){10}} \qbezier(0,25)(2,27)(4,29) \qbezier(6,31)(8,33)(10,35)
\put(0,75){\line(0,1){6}}
\put(10,75){\line(0,1){4}}
\qbezier(10,35)(9,36)(7,38) \qbezier(0,35)(1,36)(3,38)
\qbezier[5](5,40)(5,45)(5,50)
\qbezier(10,55)(9,54)(7,52) \qbezier(0,55)(1,54)(3,52)
\put(10,55){\line(-1,1){10}} \qbezier(0,55)(2,57)(4,59) \qbezier(6,61)(8,63)(10,65)
\put(10,65){\line(-1,1){10}} \qbezier(0,65)(2,67)(4,69) \qbezier(6,71)(8,73)(10,75)
\put(30,-10){\line(0,1){5}} \put(40,-10){\line(0,1){5}}
\put(30,-5){\line(1,1){10}} \qbezier(40,-5)(38,-3)(36,-1) \qbezier(34,1)(32,3)(30,5)
\put(30,5){\line(0,1){10}} \put(40,5){\line(0,1){10}}
\put(30,15){\line(1,1){10}} \qbezier(40,15)(38,17)(36,19) \qbezier(34,21)(32,23)(30,25)
\put(30,25){\line(1,1){10}} \qbezier(40,25)(38,27)(36,29) \qbezier(34,31)(32,33)(30,35)
\qbezier(30,35)(31,36)(33,38) \qbezier(40,35)(39,36)(37,38)
\qbezier[5](35,40)(35,45)(35,50)
\qbezier(30,55)(31,54)(33,52) \qbezier(40,55)(39,54)(37,52)
\put(30,55){\line(1,1){10}} \qbezier(40,55)(38,57)(36,59) \qbezier(34,61)(32,63)(30,65)
\put(30,65){\line(1,1){10}} \qbezier(40,65)(38,67)(36,69) \qbezier(34,71)(32,73)(30,75)
\put(30,75){\line(0,1){3}} \put(40,75){\line(0,1){3}}
\put(30,82){\line(0,1){3}} \put(40,82){\line(0,1){3}}
\put(30,85){\line(1,1){10}} \qbezier(40,85)(38,87)(36,89) \put(34,91){\line(-1,1){14}} 
\put(40,95){\line(0,1){5}}
\put(50,81){\line(0,1){24}}
\put(60,10){\line(0,1){5}} \put(70,-10){\line(0,1){25}}
\put(60,15){\line(1,1){10}} \qbezier(70,15)(68,17)(66,19) \qbezier(64,21)(62,23)(60,25)
\put(60,25){\line(1,1){10}} \qbezier(70,25)(68,27)(66,29) \qbezier(64,31)(62,33)(60,35)
\qbezier(60,35)(61,36)(63,38) \qbezier(70,35)(69,36)(67,38)
\qbezier[5](65,40)(65,45)(65,50)
\qbezier(60,55)(61,54)(63,52) \qbezier(70,55)(69,54)(67,52)
\put(60,55){\line(1,1){10}} \qbezier(70,55)(68,57)(66,59) \qbezier(64,61)(62,63)(60,65)
\put(60,65){\line(1,1){10}} \qbezier(70,65)(68,67)(66,69) \qbezier(64,71)(62,73)(60,75)
\put(60,75){\line(0,1){4}} \put(70,75){\line(0,1){25}}
\put(0,-10){\line(1,0){30}} \put(40,-10){\line(1,0){30}}
\put(10,10){\line(1,0){17}} \put(33,10){\line(1,0){4}} \put(43,10){\line(1,0){17}}
\mythicklines
\qbezier(10,65)(20,65)(40,65)
\qbezier(40,65)(55,60)(40,55) \qbezier(40,55)(45,55)(47,52)
\qbezier[3](47,40)(47,45)(47,50)
\qbezier(40,35)(45,35)(47,38) \qbezier(40,25)(55,30)(40,35)
\qbezier(40,25)(50,25)(70,25)
\qbezier(70,25)(85,30)(70,35) \qbezier(70,35)(75,35)(77,38)
\qbezier[3](77,40)(77,45)(77,50)
\qbezier(70,55)(85,60)(70,65) \qbezier(70,55)(75,55)(77,52)
\qbezier(70,65)(100,115)(50,115) \qbezier(50,115)(0,115)(30,85)
\qbezier(30,85)(-10,105)(-15,45)
\qbezier(-15,45)(-10,-15)(30,5)
\qbezier(30,5)(22,5)(20,10) \qbezier(10,25)(18,25)(20,10)
\qbezier(10,35)(15,35)(17,38) \qbezier(10,25)(25,30)(10,35)
\qbezier(10,65)(25,60)(10,55) \qbezier(10,55)(15,55)(17,52)
\qbezier[3](17,40)(17,45)(17,50)
\put(-15,92){(d)}
\thicklines
\qbezier[4](-3,8)(0,8)(3,8)\qbezier[4](-3,2)(0,2)(3,2)\qbezier[4](-3,2)(-3,5)(-3,8)\qbezier[4](3,2)(3,5)(3,8)
\qbezier[4](7,68)(10,68)(13,68)\qbezier[4](7,22)(10,22)(13,22)\qbezier[30](7,22)(7,45)(7,68)\qbezier[30](13,22)(13,45)(13,68)
\qbezier[4](17,108)(20,108)(23,108)\qbezier[4](17,8)(20,8)(23,8)\qbezier[66](17,8)(17,58)(17,108)\qbezier[66](23,8)(23,58)(23,108)
\qbezier[4](27,88)(30,88)(33,88)\qbezier[4](27,2)(30,2)(33,2)\qbezier[58](27,2)(27,45)(27,88)\qbezier[58](33,2)(33,45)(33,88)
\qbezier[4](37,68)(40,68)(43,68)\qbezier[4](37,22)(40,22)(43,22)\qbezier[30](37,22)(37,45)(37,68)\qbezier[30](43,22)(43,45)(43,68)
\qbezier[4](57,28)(60,28)(63,28)\qbezier[4](57,22)(60,22)(63,22)\qbezier[4](57,22)(57,25)(57,28)\qbezier[4](63,22)(63,25)(63,28)
\qbezier[4](67,68)(70,68)(73,68)\qbezier[4](67,22)(70,22)(73,22)\qbezier[30](67,22)(67,45)(67,68)\qbezier[30](73,22)(73,45)(73,68)
\end{picture}
\vspace{-3mm}
\caption{Arc presentations of $P(-p,q,r)$ with $p\geq3$ and  $4\le q\le r$}\label{fig:(-p,q,r)}
\end{figure}
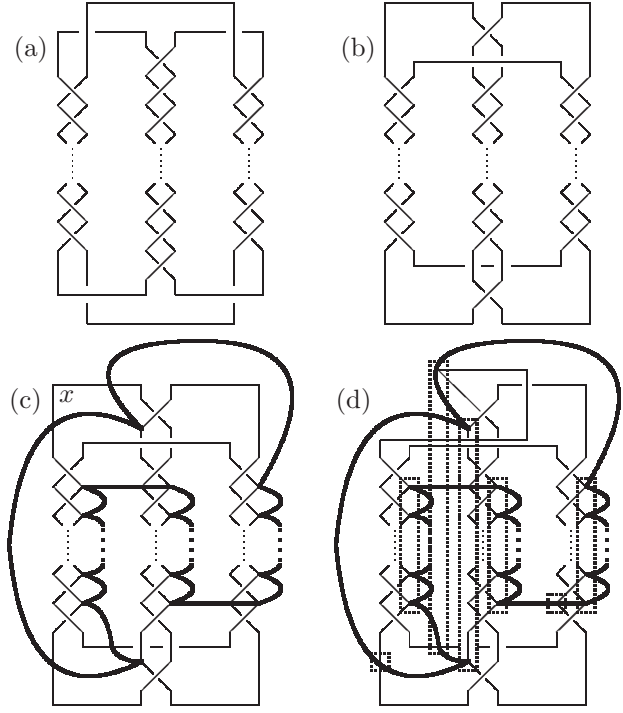

\section{The Kauffman polynomial of the pretzel knots $P(-p,q,r)$}
\newcommand{\ANG}[1]{\langle #1\rangle}
\newcommand{\SQU}[1]{\left[#1\right]}

For any link diagram $D$, the polynomial $\Lambda_D$ is of the form
$$\Lambda_D(a,z)=\sum_{i=m}^n f_i(z)a^i$$
where $m,n$ are integers with $m\le n$, and $f_i(z)$'s are polynomials in $z$ with integer coefficients such that $f_m(z)\ne0$ and $f_n(z)\ne0$.
To simplify our computation of $\operatorname{spread}_a(\Lambda_D)$ we use the notations
$$f_m(z)=\langle k_m z^{h_m}\rangle$$
$$f_n(z)=\langle k_n z^{h_n}\rangle$$
$$\sum_{i=m}^n f_i(z)a^i=\SQU{\ANG{k_n z^{h_n}}a^{n},\ANG{k_m z^{h_m}}a^{m}},\quad(m<n)$$
where $k_mz^{h_m}$ and $k_nz^{h_n}$ are the highest degree terms in $f_m(z)$ and $f_n(z)$, respectively. For example, we write
$$z(z^2-1)a^{-1}+z^2a^{-2}-2za^{-3}=\SQU{\ANG{z^3} a^{-1},\ANG{-2z} a^{-3}}.$$
We also use the notation $\Lambda_{(p_1,p_2,\ldots,p_n)}$ for $\Lambda_D$ when $D=P(p_1,p_2,\ldots,p_n)$.
	
\begin{figure}[ht]
\centering
\setlength{\unitlength}{1.0pt}\small
\mythinlines
\begin{picture}(50,80)
\put(5,80){\line(1,0){40}}
\qbezier(5,80)(0,80)(0,75)\qbezier(5,70)(0,70)(0,75)
\qbezier(45,80)(50,80)(50,75)\qbezier(45,70)(50,70)(50,75)
\put(15,70){\line(1,0){20}}
\put( 5,10){\line(0,1){5}} \put(15,10){\line(0,1){5}}
\put(15,15){\line(-1,1){10}} \qbezier(5,15)(7,17)(8,18) \qbezier(12,22)(13,23)(15,25)
\put(15,25){\line(-1,1){10}} \qbezier(5,25)(7,27)(8,28) \qbezier(12,32)(13,33)(15,35)
\put( 9,35){$\vdots$}\put(12,39){$-m$}
\put(15,45){\line(-1,1){10}} \qbezier(5,45)(7,47)(8,48) \qbezier(12,52)(13,53)(15,55)
\put(15,55){\line(-1,1){10}} \qbezier(5,55)(7,57)(8,58) \qbezier(12,62)(13,63)(15,65)
\put( 5,65){\line(0,1){5}} \put(15,65){\line(0,1){5}}
\put(35,10){\line(0,1){60}} \put(45,10){\line(0,1){60}}
\put(5,0){\line(1,0){40}}
\qbezier(5,0)(0,0)(0,5)\qbezier(5,10)(0,10)(0,5)
\qbezier(45,0)(50,0)(50,5)\qbezier(45,10)(50,10)(50,5)
\put(15,10){\line(1,0){20}}
\end{picture}
\qquad
\begin{picture}(50,80)
\put(5,80){\line(1,0){40}}
\qbezier(5,80)(0,80)(0,75)\qbezier(5,70)(0,70)(0,75)
\qbezier(45,80)(50,80)(50,75)\qbezier(45,70)(50,70)(50,75)
\put(15,70){\line(1,0){20}}
\put( 5,10){\line(0,1){60}} \put(15,10){\line(0,1){60}}
\put(35,10){\line(0,1){5}} \put(45,10){\line(0,1){5}}
\put(35,15){\line(1,1){10}} \qbezier(45,15)(43,17)(42,18) \qbezier(38,22)(37,23)(35,25)
\put(35,25){\line(1,1){10}} \qbezier(45,25)(43,27)(42,28) \qbezier(38,32)(37,33)(35,35)
\put(39,35){$\vdots$}\put(47,39){$n$}
\put(35,45){\line(1,1){10}} \qbezier(45,45)(43,47)(42,48) \qbezier(38,52)(37,53)(35,55)
\put(35,55){\line(1,1){10}} \qbezier(45,55)(43,57)(42,58) \qbezier(38,62)(37,63)(35,65)
\put(35,65){\line(0,1){5}} \put(45,65){\line(0,1){5}}
\put(5,0){\line(1,0){40}}
\qbezier(5,0)(0,0)(0,5)\qbezier(5,10)(0,10)(0,5)
\qbezier(45,0)(50,0)(50,5)\qbezier(45,10)(50,10)(50,5)
\put(15,10){\line(1,0){20}}
\end{picture}
\qquad
\begin{picture}(50,80)
\put(5,80){\line(1,0){40}}
\qbezier(5,80)(0,80)(0,75)\qbezier(5,70)(0,70)(0,75)
\qbezier(45,80)(50,80)(50,75)\qbezier(45,70)(50,70)(50,75)
\put(15,70){\line(1,0){20}}
\put( 5,10){\line(0,1){5}} \put(15,10){\line(0,1){5}}
\put(15,15){\line(-1,1){10}} \qbezier(5,15)(7,17)(8,18) \qbezier(12,22)(13,23)(15,25)
\put(15,25){\line(-1,1){10}} \qbezier(5,25)(7,27)(8,28) \qbezier(12,32)(13,33)(15,35)
\put( 9,35){$\vdots$}\put(12,39){$-m$}
\put(15,45){\line(-1,1){10}} \qbezier(5,45)(7,47)(8,48) \qbezier(12,52)(13,53)(15,55)
\put(15,55){\line(-1,1){10}} \qbezier(5,55)(7,57)(8,58) \qbezier(12,62)(13,63)(15,65)
\put( 5,65){\line(0,1){5}} \put(15,65){\line(0,1){5}}
\put(35,10){\line(0,1){5}} \put(45,10){\line(0,1){5}}
\put(35,15){\line(1,1){10}} \qbezier(45,15)(43,17)(42,18) \qbezier(38,22)(37,23)(35,25)
\put(35,25){\line(1,1){10}} \qbezier(45,25)(43,27)(42,28) \qbezier(38,32)(37,33)(35,35)
\put(39,35){$\vdots$}\put(47,39){$n$}
\put(35,45){\line(1,1){10}} \qbezier(45,45)(43,47)(42,48) \qbezier(38,52)(37,53)(35,55)
\put(35,55){\line(1,1){10}} \qbezier(45,55)(43,57)(42,58) \qbezier(38,62)(37,63)(35,65)
\put(35,65){\line(0,1){5}} \put(45,65){\line(0,1){5}}
\put(5,0){\line(1,0){40}}
\qbezier(5,0)(0,0)(0,5)\qbezier(5,10)(0,10)(0,5)
\qbezier(45,0)(50,0)(50,5)\qbezier(45,10)(50,10)(50,5)
\put(15,10){\line(1,0){20}}
\end{picture}
\caption{Links $P(-m,0)$, $P(0,n)$, and $P(-m,n)$ }\label{fig:P(-m,n)}
\end{figure}
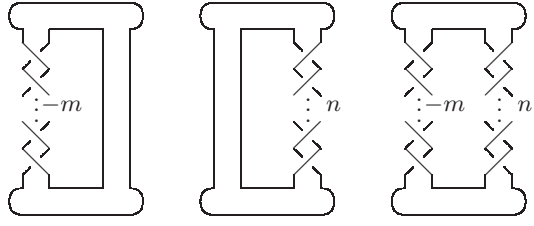

Let the polynomial $\Phi_i(z)$ be defined by $\Phi_0(z)=1, \Phi_1(z)=z$ and
$$\Phi_{i+1}(z)=z\Phi_i(z)-\Phi_{i-1}(z).\footnote{$\Phi_i(2z)$'s are the Chebyshev's polynomials of the second type.}$$

\begin{lemma}\label{lem:spread_1}
Let $m,n$ be nonnegative integers. Then
$$
\Lambda_{(-m,n)}=
\begin{cases}
\ \SQU{\ANG{z}a^{k-1},\ANG{z^{k-1}}a^{-1}}   &\text{if\/ } k=m-n>1 \\
\ a^{-1}                                                         &\text{if\/ } k=m-n=1 \\
\ z^{-1}a-1+z^{-1}a^{-1}                                         &\text{if\/ } k=n-m=0   \\
\ a                                                              &\text{if\/ } k=n-m=1 \\
\ \SQU{\ANG{z^{k-1}} a,\ANG{z} a^{-(k-1)}}     &\text{if\/ } k=n-m>1 \\
\end{cases}
$$
More precisely,
$$
\Lambda_{(-m,n)}=
\begin{cases}
\ \SQU{za^{m-1},z^{-1}\Phi_m(z)a^{-1}}   &\text{if\/ } m\geq 3,\ \ n=0 \\
\ \SQU{z^{-1}\Phi_n(z)a,za^{-(n-1)}}     &\text{if\/ } m=0,\ \ n\geq 3 \\
\end{cases}
$$
\end{lemma}

\begin{proof}
Using the skein relations K1, K2, and K3, we have the followings
$$\Lambda_{(-1,0)}=a^{-1},\ \  \Lambda_{(0,1)}=a$$
$$\Lambda_{(0,0)}=z^{-1}a-1+z^{-1}a^{-1}.$$

Let $D_v$ be a vertical integer tangle which has $|v|$ times half-twists in the positive or negative direction according to the sign of $v$ for an integer $v$. $D_{+1}$, $D_{-1}$, $D_{0}$, and $D_{\infty}$ are exemplified in ${\rm (K2)}$.  In the sublemma below, we use the following simplified notations:
$$\Lambda_v=\Lambda_{D_v},\ \ \Lambda_{\infty}=\Lambda_{D_{\infty}}.$$

\begin{sublem} For $m\ge2$ and $n\ge 2$, we have
\begin{eqnarray}
\Lambda_{-m} &=& -\Phi_{m-2}(z)\,\Lambda_0+\Phi_{m-1}(z)\,\Lambda_{-1}+\left(\sum_{i=0}^{m-2}z\Phi_i(z)a^{m-1-i}\right)\,\Lambda_\infty \label{eqn:-m} \\
\Lambda_{n} &=& -\Phi_{n-2}(z)\,\Lambda_0+\Phi_{n-1}(z)\,\Lambda_{+1}+\left(\sum_{j=0}^{n-2}z\Phi_j(z)a^{-(n-1-j)}\right)\,\Lambda_\infty
\label{eqn:n}
\end{eqnarray}
\end{sublem}

\begin{pfs} As shown below, the equation \emph{(\ref{eqn:-m})} holds when $m=2, 3$.
$$
\begin{aligned}
\Lambda_{-2} &= -\Lambda_0+z\Lambda_{-1}+za\Lambda_\infty \\
&= -\Phi_{0}(z)\,\Lambda_0+\Phi_{1}(z)\,\Lambda_{-1}+z\Phi_0(z)a\,\Lambda_\infty\\
\Lambda_{-3} &= -\Lambda_{-1}+z\Lambda_{-2}+za^2\Lambda_\infty \\
&= -\Phi_0(z)\,\Lambda_{-1}+z\left(-\Phi_{0}(z)\,\Lambda_0+\Phi_{1}(z)\,\Lambda_{-1}+z\Phi_0(z)a\,\Lambda_\infty\right)+za^2\Lambda_\infty \\
&= -\Phi_{1}(z)\,\Lambda_0+\Phi_{2}(z)\,\Lambda_{-1}+\left(\sum_{i=0}^1z\Phi_i(z)a^{2-i}\right)\,\Lambda_\infty\\
\end{aligned}
$$

Suppose that the equation \emph{(\ref{eqn:-m})} holds for $2\le m\le k-1$ for some $k\ge4$. Then
$$
\begin{aligned}
\Lambda_{-k} &= -\Lambda_{-(k-2)}+z\Lambda_{-(k-1)}+za^{k-1}\Lambda_\infty \\
&= -\left(-\Phi_{k-4}(z)\,\Lambda_0+\Phi_{k-3}(z)\,\Lambda_{-1}+\left(\sum_{i=0}^{k-4}z\Phi_i(z)a^{k-3-i}\right)\,\Lambda_\infty\right)\\
& \phantom{={}}+z\left(-\Phi_{k-3}(z)\,\Lambda_0+\Phi_{k-2}(z)\,\Lambda_{-1}+\left(\sum_{i=0}^{k-3}z\Phi_i(z)a^{k-2-i}\right)\,\Lambda_\infty\right) +za^{k-1}\Lambda_\infty \\
&= -\Phi_{k-2}(z)\,\Lambda_0+\Phi_{k-1}(z)\,\Lambda_{-1}+\left(\sum_{i=0}^{k-2}z\Phi_i(z)a^{k-1-i}\right)\,\Lambda_\infty\\
\end{aligned}
$$
This proves that the equation \emph{(\ref{eqn:-m})} holds for $m\ge2$.
In a similar manner, we can prove that the equation \emph{(\ref{eqn:n})} holds for $n\ge2$. \qed
\end{pfs}

For $m\geq 2$ and $n\geq 2$, using the equations (\ref{eqn:-m}) and (\ref{eqn:n}) on $m$ and $n$ respectively, we obtain
$$
\begin{aligned}
 \Lambda_{(-m,0)}
 &= -\Phi_{m-2}(z)\Lambda_{(0,0)}+\Phi_{m-1}(z)\Lambda_{(-1,0)}+\sum_{i=0}^{m-2}z\Phi_{i}(z)a^{m-1-i} \\
 &=\begin{cases}
                  \SQU{(z-z^{-1})a,(z-z^{-1})a^{-1}}                                                                                 &\text{if\/ } m=2\\
                  \SQU{za^{m-1},\{\Phi_{m-1}(z)-z^{-1}\Phi_{m-2}(z)\}a^{-1}}
                  &\text{if\/ } m>2
  \end{cases}\\
 &=\begin{cases}
                  \SQU{z^{-1}\Phi_2(z)a,z^{-1}\Phi_2(z)a^{-1}}                                                                                 &\text{if\/ } m=2\\
                  \SQU{za^{m-1},z^{-1}\Phi_m(z)a^{-1}}
                  &\text{if\/ } m>2
  \end{cases}
\end{aligned}
$$
$$
\begin{aligned}
 \Lambda_{(0,n)}
 &= -\Phi_{n-2}(z)\Lambda_{(0,0)}+\Phi_{n-1}(z)\Lambda_{(0,1)}+\sum_{j=0}^{n-2}z\Phi_{j}(z)a^{-(n-1-j)} \\
 &=\begin{cases}
                  \SQU{(z-z^{-1})a,(z-z^{-1})a^{-1}}                                                                                 &\text{if\/ } n=2\\
                  \SQU{\{\Phi_{n-1}(z)-z^{-1}\Phi_{n-2}(z)\}a,za^{-(n-1)}}
                  &\text{if\/ } n>2
  \end{cases}\\
 &=\begin{cases}
                  \SQU{z^{-1}\Phi_2(z)a,z^{-1}\Phi_2(z)a^{-1}}                                                                                 &\text{if\/ } n=2\\
                  \SQU{z^{-1}\Phi_n(z)a,za^{-(n-1)}}
                  &\text{if\/ } n>2
  \end{cases}
\end{aligned}
$$
Since $\Lambda_D(a,z)$ is an invariant under regular isotopy of diagrams, we have
$$
\Lambda_{(-m,n)}=
\begin{cases}
\ \Lambda_{(-k,0)}     &\text{if\/ } k=m-n\ge 1 \\
\ \Lambda_{(0,0)}      &\text{if\/ } k=m-n=0   \\
\ \Lambda_{(0,k)}      &\text{if\/ } k=n-m\ge 1 \\
\end{cases}
$$
This completes the proof.
\end{proof}

\begin{lemma}\label{lem:spread_2}
$$
\begin{aligned}
\Lambda_{(-p,0,r)}&=
\begin{cases}
\ 1                                                     &\text{if\/ } p=r=1 \\
\ \SQU{\ANG{z} a^p,\ANG{z^{p-1}} a^0}     &\text{if\/ } p>1,\ r=1   \\
\ \SQU{\ANG{z^r} a^p,\ANG{z^{p}} a^{-r}}  &\text{if\/ } p>1,\ r>1 \\
\end{cases}\\
\Lambda_{(-p,1,r)}&=
\begin{cases}
\ a^2      &\text{if\/ } p=0,\ r=1 \\
\ \SQU{\ANG{z^{r-1}} a^2,\ANG{z} a^{-(r-2)}}      &\text{if\/ } p=0,\ r>1   \\
\ a^{-r}     &\text{if\/ } p=1,\ r\ge1 \\
\ \SQU{\ANG{z^{r+1}} a^p,\ANG{z^{p-1}} a^{-r}}     &\text{if\/ } p>1,\ r\ge1 \\
\end{cases}
\end{aligned}
$$
More precisely, for $p\geq 3$ and $r\geq 3$
$$\Lambda_{(-p,0,r)}=\SQU{\Phi_r(z)a^{p},\ \Phi_p(z)a^{-r}} $$
and for $r\geq 2$
$$
\begin{aligned}
\Lambda_{(-p,1,r)}&=\begin{cases}
                  \SQU{z^{-1}\Phi_{r+2}(z)a^2,\ \Phi_{1}(z)a^{-r}}                                                                                   &\text{if\/ } p=2\\
                  \SQU{\Phi_{r+1}(z)a^{p},\ \Phi_{p-1}(z)a^{-r}}
                  &\text{if\/ } p>2
     \end{cases}
\end{aligned}
$$
\end{lemma}

\begin{proof}
Since $P(-p,0,r)=P(-p,0)\sharp P(0,r)$, we have $\Lambda_{(-p,0,r)}=\Lambda_{(-p,0)}\,\Lambda_{(0,r)}$.
Therefore the formula about $\Lambda_{(-p,0,r)}$ follows from Lemma~\ref{lem:spread_1}.

Now we consider the formula about about $\Lambda_{(-p,1,r)}$. Three cases with $p=0$ or $ p=1$ follow from K1, K2, K3, and Lemma~\ref{lem:spread_1}. The other case is derived by the equation {(\ref{eqn:-m})}.
For $p\geq 2$, we have
$$
\begin{aligned}
 \Lambda_{(-p,1,r)}
 &= -\Phi_{p-2}(z) \Lambda_{(0,1,r)}+\Phi_{p-1}(z)\Lambda_{(-1,1,r)}+\sum_{i=0}^{p-2}z\Phi_{i}(z)a^{p-1-i}\Lambda_{(1,r)}\\
 &= -\Phi_{p-2}(z)a \Lambda_{(0,r)}+\Phi_{p-1}(z)a^{-r}+\sum_{i=0}^{p-2}z\Phi_{i}(z)a^{p-1-i}\Lambda_{(0,r+1)}\\
 &=\begin{cases}
                  \SQU{\{-z^{-1}\Phi_{r}(z)+\Phi_{r+1}(z)\}a^2,\ \Phi_{1}(z)a^{-r}}                                                                                   &\text{if\/ } p=2\\
                  \SQU{\Phi_{r+1}(z)a^{p},\ \Phi_{p-1}(z)a^{-r}}
                  &\text{if\/ } p>2
     \end{cases}\\
 &=\begin{cases}
                  \SQU{z^{-1}\Phi_{r+2}(z)a^2,\ \Phi_{p-1}(z)a^{-r}}                                                                                   &\text{if\/ } p=2\\
                  \SQU{\Phi_{r+1}(z)a^{p},\ \Phi_{p-1}(z)a^{-r}}
                  &\text{if\/ } p>2
     \end{cases}
\end{aligned}
$$
This completes the proof.
\end{proof}

\begin{proposition}\label{prop:spr(-p2r)=p+r}
$\operatorname{spread}_a (\Lambda_{(-p, 2, r)}(a,z)) =p+r$ for $p\geq3$, and\/ $r\geq3$.
\end{proposition}
	
\begin{proof}
For $p\ge3$, we show that
\begin{equation}\label{eqn:lambda(-p,2,r)}
\Lambda_{(-p,2,r)}=
\begin{cases}
\SQU{\ANG{z^{3}}a^p,\ANG{z^{p-1}}a^{-2}} & (\text{if\/ } r=1)\\
\SQU{\ANG{z^{4}}a^p,\ANG{3z^{p-2}}a^{-2}} & (\text{if\/ } r=2)\\
\SQU{\ANG{z^{r+2}}a^p,\ANG{z^{p-2}}a^{-r}} & (\text{if\/ } r\ge3)
\end{cases}
\end{equation}

Using K1, K2, K3 and Lemmas~\ref{lem:spread_1} and \ref{lem:spread_2}, we obtain
\begin{align*}
\Lambda_{(-p,2,1)}&=-\Lambda_{(-p,0, 1)}+z\Lambda_{(-p,1, 1)}+za^{-1}\Lambda_{(-p,1)}\\
                  &=-\SQU{\ANG{z}a^{p},\ANG{z^{p-1}}a^{0}} + z\SQU{\ANG{z^2}a^p,\ANG{z^{p-1}}a^{-1}}
                         +za^{-1}\SQU{\ANG{z}a^{p-2},\ANG{z^{p-2}}a^{-1}}\\
                  &=\SQU{\ANG{z^3}a^p,\ANG{z^{p-1}}a^{-2}},\\
\Lambda_{(-p,2,2)}&=-\Lambda_{(-p,0, 2)}+z\Lambda_{(-p,1, 2)}+za^{-1}\Lambda_{(-p,2)}\\
                  &=-\SQU{za^{p-1}, z^{-1}\Phi_p(z)a^{-1}}\SQU{z^{-1}\Phi_2(z)a,z^{-1}\Phi_2(z)a^{-1}}\\
                  &    \qquad\qquad   +       z\SQU{\Phi_3(z)a^p,\Phi_{p-1}(z)a^{-2}}\\
                  &    \qquad\qquad   +  za^{-1} \begin{cases}
                                              a^{-1} & (\text{if\/ } p=3)\\
                                              \SQU{z^{-1}\Phi_2(z)a,z^{-1}\Phi_2(z)a^{-1}} & (\text{if\/ } p=4)\\
                                              \SQU{za^{p-3}, z^{-1}\Phi_{p-2}(z)a^{-1}} & (\text{if\/ } p\ge5)\\
                                              \end{cases}\\
                  &=-\SQU{\Phi_2(z)a^p,z^{-2}\Phi_2(z)\Phi_p(z)a^{-2}}+z\SQU{\Phi_3(z)a^p,\Phi_{p-1}(z)a^{-2}}\\
                  &    \qquad\qquad   +       \begin{cases}
                                              za^{-2} & (\text{if\/ } p=3)\\
                                              za^{-1}\SQU{z^{-1}\Phi_2(z)a,z^{-1}\Phi_2(z)a^{-1}} & (\text{if\/ } p=4)\\
                                              za^{-1}\SQU{za^{p-3}, z^{-1}\Phi_{p-2}(z)a^{-1}} & (\text{if\/ } p\ge5)\\
                                              \end{cases}\\
                  &=\SQU{\Phi_4(z)a^p,\ANG{3z^{p-2}}a^{-2}},
\end{align*}
which prove the first two cases of (\ref{eqn:lambda(-p,2,r)}). Now we show the third case of (\ref{eqn:lambda(-p,2,r)}) by an induction on $r$.
For $r=3$, we have
\begin{align*}
\Lambda_{(-p,2,3)}&=-\Lambda_{(-p,2, 1)}+z\Lambda_{(-p,2, 2)}+za^{-2}\Lambda_{(-p,2)}\\
                  &=-\SQU{\ANG{z^3}a^p,\ANG{z^{p-1}}a^{-2}}+z\SQU{\ANG{z^4}a^p,\ANG{3z^{p-2}}a^{-2}}\\
                  &    \qquad\qquad   +       \begin{cases}
                                              za^{-3} & (\text{if\/ } p=3)\\
                                              \SQU{\ANG{z^2}a^{p-5},\ANG{z^{p-2}}a^{-3}} & (\text{if\/ } p\ge4)\\
                                              \end{cases}\\
                  &=\SQU{\ANG{z^5}a^p,\ANG{z^{p-2}}a^{-3}},
\end{align*}
and for $r\ge4$, inductively, we have
\begin{align*}
\Lambda_{(-p,2,r)}&=-\Lambda_{(-p,2,r-2)}+z\Lambda_{(-p,2,r-1)}+za^{-(r-1)}\Lambda_{(-p,2)}\\
                  &=-\SQU{\ANG{z^r}a^p,\ANG{*}a^{-(r-2)}}+z\SQU{\ANG{z^{r+1}}a^p,\ANG{z^{p-2}}a^{-(r-1)}} \\
                  &    \qquad\qquad   +       \begin{cases}
                                              za^{-r} & (\text{if\/ } p=3)\\
                                              \SQU{\ANG{z^2}a^{p-r-2},\ANG{z^{p-2}}a^{-r}} & (\text{if\/ } p\ge4)\\
                                              \end{cases}\\
                  &=\SQU{\ANG{z^{r+2}}a^p,\ANG{z^{p-2}}a^{-r}}
\end{align*}
where $\ANG{*}a^{-(r-2)}$ indicates that the lowest $a$-degree of $\Lambda_{(-p,2,r-2)}$ is not smaller than $-(r-2)$.

This completes the proof.
\end{proof}
	
\begin{proposition}\label{prop:spr(-p3r)=p+r}
$\operatorname{spread}_a (\Lambda_{(-p, 3, r)}(a,z)) =p+r$ for $p\geq3$ and\/  $r\geq3$.
\end{proposition}
	
\begin{proof}
We show that
\begin{equation*}
\Lambda_{(-p,3,r)}=
\begin{cases}
\SQU{\ANG{z^6}a^p,\ANG{2z^{p-3}}a^{-3}} &(\text{if\/ } r=3)\\
\SQU{\ANG{z^{r+3}}a^p,\ANG{z^{p-3}}a^{-r}} &(\text{if\/ } r\ge4)\\
\end{cases}
\end{equation*}
Using the equation {(\ref{eqn:n})} and Lemmas~\ref{lem:spread_1} and \ref{lem:spread_2}, we obtain
$$
\begin{aligned}
\Lambda_{(-p,3,3)}&=-\Phi_1(z)\Lambda_{(-p,0,3)}+\Phi_2(z)\Lambda_{(-p,1,3)}+(za^{-2}+z^2a^{-1})\Lambda_{(-p,3)}\\
                  &=-\Phi_1(z)\SQU{za^{p-1}, z^{-1}\Phi_p(z)a^{-1}}\SQU{z^{-1}\Phi_3(z)a,za^{-2}}\\
                  &\qquad\qquad +\Phi_2(z)\SQU{\Phi_4a^{p},\Phi_{p-1}(z)a^{-3}}\\
                  &\qquad\qquad +(za^{-2}+z^2a^{-1})
                           \begin{cases}
                           (z^{-1}a-1+z^{-1}a^{-1}) &(\text{if\/ } p=3)\\
                           a^{-1} &(\text{if\/ } p=4)\\
                           \SQU{\ANG{z}a^{p-4},\ANG{z^{p-4}}a^{-1}} &(\text{if\/ } p\ge5)\\
                           \end{cases}\\
                  &=\SQU{\ANG{z^6}a^p,\Phi_{p-3}(z)a^{-3}}+
                           \begin{cases}
                           \SQU{z,a^{-3}} &(\text{if\/ } p=3)\\
                           \SQU{z^2a^{-2},za^{-3}} &(\text{if\/ } p=4)\\
                           \SQU{\ANG{z^3}a^{p-5},\ANG{z^{p-3}}a^{-3}} &(\text{if\/ } p\ge5)\\
                           \end{cases}\\         
                  &=\SQU{\ANG{z^{6}}a^p,\ANG{2z^{p-3}}a^{-3}}.\\
\end{aligned}
$$
Using K1, K2, K3 and the results above, we obtain
$$
\begin{aligned}
\Lambda_{(-p,3,4)}&=-\Lambda_{(-p,2, 3)}+z\Lambda_{(-p,3, 3)}+za^{-3}\Lambda_{(-p,3)} \qquad(\because \Lambda_{(-p,3,2)}=\Lambda_{(-p,2, 3)})\\
                  &=-\SQU{\ANG{z^5}a^p,\ANG{z^{p-2}}a^{-2}}+z\SQU{\ANG{z^6}a^p,\ANG{2z^{p-3}}a^{-3}}\\
                  &\qquad\qquad+za^{-3}
                           \begin{cases}
                           (z^{-1}a-1+z^{-1}a^{-1}) &(\text{if\/ } p=3)\\
                           a^{-1} &(\text{if\/ } p=4)\\
                           \SQU{\ANG{z}a^{p-4},\ANG{z^{p-4}}a^{-1}} &(\text{if\/ } p\ge5)\\
                           \end{cases}\\
                  &=\SQU{\ANG{z^7}a^{p},\ANG{z^{p-3}}a^{-4}},
\end{aligned}
$$
and, for $r\ge5$, inductively, we have
$$
\begin{aligned}
\Lambda_{(-p,3,r)}&=-\Lambda_{(-p,3, r-2)}+z\Lambda_{(-p,3, r-1)}+za^{-(r-1)}\Lambda_{(-p,3)}\\
                  &=-\SQU{\ANG{z^{r+1}}a^p,\ANG{*}a^{-(r-2)}}+z\SQU{\ANG{z^{r+2}}a^p,\ANG{z^{p-3}}a^{-(r-1)}}\\
                  &\qquad\qquad+za^{-(r-1)}
                           \begin{cases}
                           (z^{-1}a-1+z^{-1}a^{-1}) &(\text{if\/ } p=3)\\
                           a^{-1} &(\text{if\/ } p=4)\\
                           \SQU{\ANG{z}a^{p-4},\ANG{z^{p-4}}a^{-1}} &(\text{if\/ } p\ge5)\\
                           \end{cases}\\
                  &=\SQU{\ANG{z^{r+3}}a^{p},\ANG{z^{p-3}}a^{-r}}
\end{aligned}
$$
This completes the proof.
\end{proof}

\begin{proposition}\label{prop:spr(-34r)=r+1}
$\operatorname{spread}_a (\Lambda_{(-3, 4, r)}(a,z)) =r+1$ for $r\ge 7$.
\end{proposition}
	
\begin{proof}
$$
\begin{aligned}
\Lambda_{(-3,4,0)}=\Lambda_{(-3,0)}\Lambda_{(0,4)}
                  =\SQU{\Phi_{4}(z)a^3,\Phi_{3}(z)a^{-4}}\\
\end{aligned}
$$
Using the equation (\ref{eqn:-m}) and K1, K3 and Lemma~\ref{lem:spread_1}, we obtain
$$
\begin{aligned}
\Lambda_{(-3,4,1)}&=-\Phi_{1}(z)\Lambda_{(0,4,1)}+\Phi_{2}(z)\Lambda_{(-1,4,1)}+\left(\sum_{i=0}^{1}z\Phi_i(z)a^{2-i}\right)\Lambda_{(4,1)}\\
                  &=-za\Lambda_{(0,4)}+\Phi_{2}(z)a^{-4}+(za^2+z^2a)\Lambda_{(0,5)}\\
                  &=-za\SQU{z^{-1}\Phi_{4}(z)a,za^{-3}}+\Phi_{2}(z)a^{-4}+(za^2+z^2a)\SQU{z^{-1}\Phi_{5}(z)a,za^{-4}}\\
                  &=\SQU{\Phi_{5}(z)a^3,\Phi_{2}(z)a^{-4}}\\
\end{aligned}
$$
Using the equation (\ref{eqn:n}) and the above two formulas, for $r\ge 7$ we obtain
$$
\begin{aligned}
\Lambda_{(-3,4,r)}&=-\Phi_{r-2}(z)\Lambda_{(-3,4,0)}+\Phi_{r-1}(z)\Lambda_{(-3,4,1)}+\sum_{i=0}^{r-2}z\Phi_{i}(z)a^{-(r-1-i)}\Lambda_{(-3,4)}\\
                  &=-\Phi_{r-2}(z)\SQU{\Phi_{4}(z)a^3,\Phi_{3}(z)a^{-4}}\\
                  &\phantom{={}} +\Phi_{r-1}(z)\SQU{\Phi_{5}(z)a^3,\Phi_{2}(z)a^{-4}}+\sum_{i=0}^{r-2}z\Phi_{i}(z)a^{-(r-2-i)}\\
                  &=\SQU{\ANG{z^{r+4}}a^3,za^{-(r-2)}}.
\end{aligned}
$$
This completes the proof.
\end{proof}

\begin{proposition}\label{prop:spr(-p4r)=p+r}
$\operatorname{spread}_a (\Lambda_{(-p, 4, r)}(a,z)) =p+r$ for $p\geq5$ and\/ $r\ge5$.
\end{proposition}
	
\begin{proof}
Using the equation (\ref{eqn:n}) and Lemmas~\ref{lem:spread_1} and \ref{lem:spread_2}, we obtain
$$
\begin{aligned}
\Lambda_{(-p,4,r)}&=-\Phi_{2}(z)\Lambda_{(-p,0,r)}+\Phi_{3}(z)\Lambda_{(-p,1,r)}+\sum_{i=0}^{2}z\Phi_{i}(z)a^{-(3-i)}\Lambda_{(-p,r)}\\
                  &=-(z^2-1)\SQU{\Phi_{r}(z)a^{p},\ \Phi_{p}(z)a^{-r}}\\
                  &\phantom{={}} +(z^3-2z)\SQU{\Phi_{r+1}(z)a^{p},\ \Phi_{p-1}(z)a^{-r}}\\
                  &\phantom{={}} +\{za^{-3}+z^2a^{-2}+(z^3-z)a^{-1}\}\begin{cases}
                      \ \SQU{\ANG{z}a^{k-1},\ANG{z^{k-1}}a^{-1}}   &\text{if\/ } k=p-r>1 \\
                      \ a^{-1}                                     &\text{if\/ } k=p-r=1 \\
                      \ z^{-1}a-1+z^{-1}a^{-1}                     &\text{if\/ } k=p-r=0 \\
                      \ a                                          &\text{if\/ } k=r-p=1 \\
                      \ \SQU{\ANG{z^{k-1}} a,\ANG{z} a^{-(k-1)}}   &\text{if\/ } k=r-p>1 \\
                      \end{cases}\\
                  &=\SQU{\ANG{z^{r+4}}a^{p},\{-(z^2-1)\Phi_{p}(z)+(z^3-2z)\Phi_{p-1}(z)\}a^{-r}}\\
                  &=\SQU{\ANG{z^{r+4}}a^{p},\Phi_{p-4}(z)a^{-r}}.
\end{aligned}
$$
This completes the proof.
\end{proof}

\section{Proofs of main results and comments}

Theorem~\ref{thm:-2qr} is proved by Proposition~\ref{prop:-2qr}. Table~\ref{tab:-2qr} shows that the upper bound `$c(K)-1$' for the arc index in Theorem~\ref{thm:-2qr} is best possible. It also shows that the lower bound `$\operatorname{spread}_a (F_K)+2$' in Theorem~\ref{thm:lower bound of arc index} is best possible.

\begin{table}[htb]
\renewcommand\arraystretch{1.3}
\caption{Examples of Theorem~\ref{thm:-2qr}}\label{tab:-2qr}
\begin{tabular}{|c|c|c|c|c|}
\hline
Pretzel knot $K$ & DT Name\footnotemark & $\operatorname{spread}_a (F_K)+2$ & arc index & $c(K)-1$ \\\hline
$P(-2,3,3)$ & $8n3$ &  $6$ & $7$ & 7\\\hline
$P(-2,3,5)$ & $10n21$ &  $6$ & $8$ & 9\\\hline
$P(-2,3,7)$ & $12n242$ &  $9$ & $9$ & 11\\\hline
$P(-2,5,5)$ & $12n725$ &  $10$ & $10$ & 11\\\hline
\end{tabular}
\end{table}
\footnotetext{The Dowker-Thistlethwaite name. See \cite{knotinfo}.}

The proof of Theorem~\ref{thm:-p2r} is a combination of Propositions~\ref{prop:-pqr,q=2,3} and \ref{prop:spr(-p2r)=p+r}.
The proof of Theorem~\ref{thm:-p3r} is a combination of Propositions~\ref{prop:-pqr,q=2,3} and \ref{prop:spr(-p3r)=p+r}.
The proof of Theorem~\ref{thm:-p4r} is a combination of Propositions~\ref{prop:-pqr,q>3} and \ref{prop:spr(-p4r)=p+r}.
The proof of Theorem~\ref{thm:-34r} is a combination of Propositions~\ref{prop:-pqr,q>3} and \ref{prop:spr(-34r)=r+1}.
Table~\ref{tab:-34r} shows that the upper bound `$c(K)-2$' is best possible but the lower bound `$c(K)-4$' may not be best possible.

\begin{table}[htb]
\renewcommand\arraystretch{1.3}
\caption{Examples of Theorem~\ref{thm:-34r}}\label{tab:-34r}
\begin{tabular}{|c|c|c|c|c|}
\hline
Pretzel knot $K$ & DT Name & $c(K)-4$ & arc index & $c(K)-2$ \\\hline
$P(-3,4,5)$ & $12n475$ &  $8$ & $10$ & 10\\\hline
$P(-3,4,7)$ & $14n12205$ &  $10$ & $11$ & 12\\\hline
\end{tabular}
\end{table}



\bibliographystyle{amsplain}

\end{document}